\documentclass[preprint,12pt]{elsarticle}
\usepackage{amsmath}
\usepackage{amsfonts}
\usepackage{amssymb}
\usepackage{amsthm}
\usepackage{amssymb}
\usepackage{hyperref}
\usepackage{url}
\usepackage{indentfirst}
\usepackage{tikz}
\usepackage{filecontents}
\usepackage{booktabs}
\usepackage[noend]{algpseudocode}
\usepackage[margin=1.5in]{geometry}    
\usepackage{algorithm}
\usepackage{graphicx}
\usetikzlibrary{shapes,positioning}

\newcommand{\R}{\mathbb R} 

\newtheorem{proposition}{Proposition}[section]
\newtheorem{lemma}{Lemma}[section]

\newtheorem{remark}{Remark}[section]
\def\cG{\mathcal G}
\usepackage{chngcntr}
\counterwithout{remark}{section}
\counterwithout{figure}{section}
\counterwithout{table}{section}
\counterwithout{theorem}{section}
\counterwithout{lemma}{section}
\counterwithout{proposition}{section}
\journal{Annals of Operations Research}

\begin{document}
\begin{frontmatter}



\title{On Nash-solvability of $n$-person graphical games under
Markov and a-priori realizations}



\author{Vladimir Gurvich}
\ead{
vladimir.gurvich@gmail.com}
\address{RUTCOR, Rutgers University, Piscataway, New Jersey, United States} 

\author{Mariya Naumova}
\ead{mnaumova@business.rutgers.edu} 
\address{Rutgers Business School, Rutgers University, 
Piscataway, New Jersey, United States}

\address{Piscataway, New Jersey, United States}

\begin{abstract}
We consider graphical $n$-person games
with perfect information that have no Nash equilibria in pure stationary strategies. Solving these games in stationary 
mixed  strategies, we introduce
probability distributions in all non-terminal positions.
The corresponding plays can be analyzed 
under two different basic assumptions: 
the Markov and a-priori realizations.
The former one guarantees existence of a uniform best response for each player in every situation.
Nevertheless, Nash equilibrium  
may fail to exist even in mixed strategies.
The classical Nash's theorem is not applicable, 
since Markov realizations may result in   
discontinuous limit distributions and 
expected payoffs.
Although a-priori realizations 
does not share many nice properties of Markov realizations 
(for example, the existence of uniform best responses)   
but in return, Nash's theorem is applicable.
We illustrate both realizations in details
by two examples with  $2$  and  $3$  players. 
We also survey some general results 
related to Nash-solvability, 
in pure and mixed stationary strategies, 
of stochastic $n$-person games 
with perfect information 
and $n$-person graphical games among them.
\end{abstract}


\end{frontmatter}
\section{Introduction}\label{s0}
\subsection{Graphical $n$-person games with terminal payoffs}
\label{ss00}
\subsubsection*{Game structures}
Let  $G = (V,E)$ be  a finite directed graph (digraph) whose
vertices  $v \in V$   and
directed edges (arcs)  $e \in E$  are interpreted
as positions and moves of a game
of  $n$  players  $I = \{1, \ldots, n\}$.
Furthermore, let $D$  be a partition of  $V$  into  $n+2$ subsets: $D : V = V_1 \cup \ldots \cup V_n \cup V_R \cup V_T$, interpreted as follows:
\begin{itemize}
\item[] $V_T$  are terminal positions of  $G$,
from which there are no moves,
\item[] $V_R$  are positions of chance.
\item[] $V_i$ are positions controlled by the player  $i \in I$,
\end{itemize}

For each position $v \in V_i$ player $i$ chooses 
a move from position $v$, that is, an arc  $(v, v') \in E$.

For each $v \in V_R$ we fix a probability distribution $q(v)$ 
over the set of all moves $(v, v')$  from  $v$.
In other words, we define $q = \{q(v) \mid v \in V_R\}$ by assigning real numbers $q(v,v') \geq 0$ to each arc $(v,v') \in E$ such that $\sum_{v' \mid (v,v') \in E} q(v, v') = 1$ for each  $v \in V_R$. 
Probabilities  $q(v,v') = 0$  are allowed. In this case arc $(v, v')$  can be deleted from $E$  while vertices  $v$ and $v'$  remain in  $V$.

Without loss of generality  (WLOG) 
we assume that  
for each position $v \in V \setminus V_T$  the number of moves 
from  $v$  is at least $2$. 
Indeed, if  $(v,v')$ is a forced move in $v$ then 
we can contract this edge, 
that is, we delete it and merge  $v$ and  $v'$.

The initial position  $v_0 \in V$  may be fixed or not.
A quadruple $(G, D, q, v_0)$ and triplet $(G, D, q)$
will be referred to as the {\em graphical game structure}:
\textit{initialized} and \textit{not initialized}, respectively.
To simplify our notation, we replace $(G, D, q)$  by $\cG$. 

A  game structure is called
\begin{itemize}
\item[] {\em deterministic} if  $V_R = \emptyset$;
\item[] {\em almost  deterministic}  if it is initialized, 
$V_R = \{v_0\}$, 
and there is a move $(v_0, v) \in E$  to each  
$v \in V \setminus \{V_T \cup \{v_0\}\}$;
\item[] {\em play-once} if  
$|V_i| = 1$  for all  $i \in I$, that is, if each of $n$ players controls a unique position.
\end{itemize}


\subsubsection*{Initializing extensions}
Given a non-initialized game structure  $\cG$, 
let us add to its digraph  $G$  
a new position of chance  $v_0$  and 
a  move $(v_0,v)$  to every  $v \in V \setminus V_T$.
Then, let us fix an arbitrary probability distribution  $q(v_0)$  
on these new edges.
We denote the obtained initialized game structure by  
$\cG' = (G', D', q', v_0)$
and call  $\cG'$ the  {\em initializing extension} of  $\cG$.

By definition, $\cG'$ is almost deterministic 
if and only if  $\cG$ is deterministic.

\begin{remark}
Let us note that initializing game structure $\cG$ 
is obviously equivalent with introducing 
an initial probability distribution in it, 
instead of fixing an initial position. 
 However, these two approaches differ  
in the general framework of Nash-solvability. 
To construct a deterministic 
Nash equilibrium free graphical game 
is more difficult than an almost deterministic one; 
see Subsection \ref{ss05} for the definitions and  
Section \ref{s5} for more details.
\end{remark}

\subsubsection*{Two examples} 
2- and 3-person person non-initialized deterministic
play-once game structures  $\cG_2$  and $\cG_3$ are shown in Figures 1 and 2. 
Each player  $i \in I$ controls a unique  position  $v_i$
in which (s)he has two possible moves:
(\textit{f}) to follow the cycle and (\textit{t}) to terminate in  $a_i$;
for $\cG_2$  and $\cG_3$, we have 
$n = 2$  and  $n = 3$, respectively.
The initializing extensions $\cG'_2$  and $\cG'_3$
are given in the same two figures.

\tikzset{ell/.style={ellipse,draw,minimum height=0.5cm,minimum width=0.8cm,inner sep=0.25cm}}

\begin{figure}[H]
\begin{center}
\resizebox{10cm}{4cm}{
\begin{tikzpicture}
\node[ell] (e1)at (-10,0) {$v_1$};
\node[ell] (e2)at (-7,0) {$v_2$};
\node[ell] (ea1)at (-10,2) {$a_1$};
\node[ell] (ea2)at (-7,2) {$a_2$};
\node[] (through2) at (-8.5, -0.5){};
\node (node1) at (-11,3) {$\cG_2:$};

\node[ell] (_e1)at (-2.5,0) {$v_1$};
\node[ell] (_e2)at (0.5,0) {$v_2$};
\node[ell] (_ea1)at (-2.5,2) {$a_1$};
\node[ell] (_ea2)at (0.5,2) {$a_2$};
\node[ell, label=below: ] (_c)at (-1,0.8) {$v_0$};
\node[] (_through2) at (-1, -0.5) {}; 
\node (_node1) at (-3.5,3) {$\cG_2':$};

\draw [->] (_c) to node[above]{\footnotesize$q_1$} (_e1);
\draw [->] (_c) to node[above]{\footnotesize$q_2$} (_e2);
\draw [->] (e1) to (e2);
\draw [->] (e1) to (ea1);
\draw [->] (e2) to (ea2);
\draw[->,rounded corners] (e2.south) |- (through2.south) -| (e1.south);

\draw [->] (_e1) to (_e2);
\draw [->] (_e1) to (_ea1);
\draw [->] (_e2) to (_ea2);
\draw[->,rounded corners] (_e2.south) |- (_through2.south) -| (_e1.south);
\end{tikzpicture}}
\end{center} 
\caption{
Game structures  $\cG_2$  and   $\cG_2'$.}
\label{Fig1}
\end{figure}

\begin{figure}[H]
\begin{center}
\resizebox{10cm}{5cm}{
\begin{tikzpicture}
\node[ell, label=above: ] (ea3)at (-8.5,4) {$a_3$};
\node[ell, label=above: ] (e3)at (-8.5,2.25) {$v_3$};
\node[ell, label=below: ] (e1)at (-10,0) {$v_1$};
\node[ell, label=below: ] (e2)at (-7,0) {$v_2$};
\node[ell, label=below: ] (ea1)at (-10,2.25) {$a_1$};
\node[ell, label=below: ] (ea2)at (-7,2.25) {$a_2$};
\node (node1) at (-11,4) {$\cG_3:$};

\node[ell, label=above: ] (_ea3)at (-1,4) {$a_3$};
\node[ell, label=above: ] (_e3)at (-1,2.25) {$v_3$};
\node[ell, label=below: ] (_e1)at (-2.5,0) {$v_1$};
\node[ell, label=below: ] (_e2)at (0.5,0) {$v_2$};
\node[ell, label=below: ] (_ea1)at (-2.5,2.25) {$a_1$};
\node[ell, label=below: ] (_ea2)at (0.5,2.25) {$a_2$};
\node[ell, label=below: ] (_c)at (-1,0.8) {$v_0$};
\node (node2) at (-3.5,4) {$\cG_3':$};

\draw [->] (e1) to  (e2);
\draw [->] (e2) to (e3);
\draw [->] (e3) to  (e1);
\draw [->] (e1) to  (ea1);
\draw [->] (e2) to  (ea2);
\draw [->] (e3) to  (ea3); 
\draw [->] (_e1) to  (_e2);
\draw [->] (_e2) to (_e3);
\draw [->] (_e3) to (_e1);
\draw [->] (_e1) to (_ea1);
\draw [->] (_e2) to (_ea2);
\draw [->] (_e3) to (_ea3);
\draw [->] (_c) to node[above]{\footnotesize$q_1$} (_e1);
\draw [->] (_c) to node[above]{\footnotesize$q_2$} (_e2);
\draw [->] (_c) to node[right]{\footnotesize$q_3$} (_e3);
\end{tikzpicture}}
\end{center} 
\caption{
Game structures  $\cG_3$  and   $\cG_3'$.}
\label{Fig2}
\end{figure}

\subsubsection*{Plays, outcomes, payoffs, and games}
A {\em play} is as a directed walk in $G$
that begins in some position  $v \in V$.
In the initialized case we assume that 
$v = v_0$, while in the non-initialized case, 
$v$  can be any position in $V \setminus V_T$.
A play is finite if and only if it ends in $V_T$.
In this case it is called  
{\em terminal},  
otherwise, it is called {\em infinite}.

Every terminal $v \in V_T$  is an {\em outcome}; 
all terminal plays ending in  $v$  
are treated as equivalent; 
they form a single outcome $v$; 
also all infinite plays are treated as equivalent; 
they form one extra outcome  $c$.
The set of outcomes will be denoted by
$A =   \{V_T \cup \{c\}\} = \{a_1, \ldots, a_m, c\}$.

A {\em payoff function} is defined as a mapping
$u : I \times A \rightarrow \R$;
the real number  $u_i(a)$  
is interpreted as the profit
of the player  $i \in I$  in case 
the outcome  $a \in A$  is realized. 
A triplet $(\cG, u, v_0)$ and a pair  $(\cG, u)$ 
will be called {\em initialized} and 
{\em non-initialized graphical games}, respectively. 
{\em Deterministic graphical} games 
were introduced in \cite{Was90} 
for the 2-person zero-sum case. 
We generalize this model allowing $n$ players and positions of chance.

\subsection{Pure stationary strategies and
normal forms of deterministic game structures}
\label{ss01}
A {\em pure stationary strategy} 
(or simply a strategy, for short) 
of a player  $i \in I$  
is a mapping $s^i : V_i \rightarrow E$  that assigns to each position
$v \in V_i$  a move  $(v, v') \in E$  from  $v$.
In other words, player  $i$  in advance makes a decision, 
how (s)he will play in each position.


An $n$-tuple  $s = \{s^i \mid i \in I\}$
of strategies of all  $n$  players is called a  {\em situation}.

If game structure $(G, D, v_0)$  is deterministic and initialized 
then each situation  $s$  uniquely defines a walk $W  = W(s)$ 
called a {\em play}. 
It begins in the initial position  $v_0$. 
Assume that  $v_0 \in V_i$.  
Then  $W$  proceeds from  $v_0$  with the move  $(v_0,v') \in E$
chosen in $v_0$  by strategy  $s^i \in s$.  
Assume that  $v' \in V_j$. 
(Equality  $i=j$  is allowed.) 
Then  $W$ proceeds from  $v'$ with the move  $(v',v'') \in E$
chosen in  $v'$  by strategy $s^j \in s$, etc.
Play  $W$  either ends in a terminal  $a  \in V_T$ 
(in which case each player  
$i \in I$  gets a profit  $u_i(a)$) 
or $W$  lasts infinitely.
Since digraph $G$  is finite, in the latter case  walk  $W$,
sooner or  later, will revisit a position, 
thus making a directed cycle.
Let us consider the first such revisiting and 
the corresponding directed cycle.
This cycle is {\em simple} 
(that is, it has no self-intersections)  
and play  $W(s)$  will repeat this cycle infinitely, 
because the players are restricted 
to their pure stationary strategies.
Such an infinite play will be called a {\em lasso}.
It consists of the initial part 
(before the first revisiting, 
which is empty if the play returns to $v_0$). 
If  $W$  is a lasso then  
each player  $i \in I$  gets a profit  $u_i(c)$, 
because we assume that all infinite plays are equivalent
and form a single outcome.

\medskip

If player  $i$  controls  $t_i$  positions
with  $\ell^i_1, \ldots, \ell^i_{t_i}$  outgoing arcs,
then  $i$  has  $k_i = \prod_{j = 1}^{t_i} \ell_j^i$  
pure stationary strategies.
In our examples $\cG_2$  and $\cG_3$ 
each player has only two such strategies.

\medskip

Given an initialized deterministic 
game structure $(G, D, v_0)$,
let  $S^i$  denote the set  of pure stationary strategies 
of player  $i \in I$ and let 
$S = S^1 \times \ldots \times S^n$  
be the direct product of  $n$  these sets.
Mapping  $g^{v_0} : S \rightarrow A$  that assigns 
to each situation  $s \in S$ 
either a terminal outcome  $v \in  V_T$, 
if play  $W(s)$  ends in  $v$, or 
the special outcome  $c$, 
if  $W(s)$  is a lasso, is called the 
{\em normal form of the 
(deterministic initialized) game structure}  $\cG$.
 
Given a non-initialized deterministic game structure $(G, D)$,
we  define its {\em normal form} as the mapping
$g : S \rightarrow 2^A$,  
where  $g(s) = \{g^v(s)  \mid  v \in V \setminus  V_T\}$.
Two examples are given in Figure \ref{Fig3}.
\begin{figure}[]
  \centering
  \begin{tikzpicture}
\node (a) at (1,2) {$t$};
\node (a) at (2.5,2) {$f$};
\node (a) at (0,0) {$f$};
\node (b) at (1,0) {$a_2\;a_2$};
\draw (0.25,-0.5) rectangle (1.75,0.5);
\node (c) at (2.5,0) {$c\;\text{ }\;c$};
\draw (1.75,-0.5) rectangle (3.25,0.5);
\node (a) at (0,1) {$t$};
\node (b) at (1,1) {$a_1\;a_2$};
\draw (0.25,0.5) rectangle (1.75,1.5);
\node (c) at (2.5,1) {$a_1\;a_1$};
\draw (1.75,0.5) rectangle (3.25,1.5);
\node (a) at (4.5,2) {$t$};
\node (a) at (5.5,2) {$f$};
\node (a) at (4,0) {$f$};
\node (b) at (4.5,0) {$\text{ }a_2$};
\draw (4.25,-0.5) rectangle (5,0.5);
\node (c) at (5.5,0) {$c\text{ }$};
\draw (5,-0.5) rectangle (5.75,0.5);
\node (a) at (4,1) {$t$};
\node (b) at (4.5,1) {$\text{ }a_1$};
\draw (4.25,0.5) rectangle (5,1.5);
\node (c) at (5.5,1) {$a_1\text{ }$};
\draw (5,0.5) rectangle (5.75,1.5);

\path[->] (4.6,0.75) edge (4.6,0.25);
\path[->] (5.4,0.25) edge (5.4,0.75);
\path[->] (4.85,0) edge (5.25,0);
\path[->] (0.75,0.75) edge (0.75,0.25);
\path[->] (2.25,0.25) edge (2.25,0.75);
\path[->] (1.25,-0.25) edge [out=-40, in=-135] (2.75,-0.25);
\path[->] (2.75,1.25) edge [out=135, in=40] (1.25,1.25);
\node (a) at (7,2) {$t$};
\node (a) at (8,2) {$f$};
\node (a) at (6.5,0) {$f$};
\node (b) at (7,0) {$\text{ }a_2$};
\draw (6.75,-0.5) rectangle (7.5,0.5);
\node (c) at (8,0) {$c\text{ }$};
\draw (7.5,-0.5) rectangle (8.25,0.5);
\node (a) at (6.5,1) {$t$};
\node (b) at (7,1) {$\text{ }a_2$};
\draw (6.75,0.5) rectangle (7.5,1.5);
\node (c) at (8,1) {$a_1\text{ }$};
\draw (7.5,0.5) rectangle (8.25,1.5);

\path[->] (7.65,1.05) edge (7.3,1.05);
\path[->] (7.9,0.25) edge (7.9,0.75);
\path[->] (7.35,0) edge (7.75,0);
\node (a) at (1.75,2.5) {$g_2$};
\node (a) at (5.15,2.5) {$g_2^{v_1}$};
\node (a) at (7.55,2.5) {$g_2^{v_2}$};
\node (a) at (-1,1.5) {$1$};
\node (b) at (-0.5,2.5) {$2$};
\path[->] (-1.25,2.25) edge (-1.25,1.5);
\path[->] (-1.25,2.25) edge (-0.5,2.25);
\node (a) at (-1.75,-1.5) {$U_2 \text{:}$};
\node (a) at (3.5,-1.5) {$u_{1}(c)>u_{1}(a_{1})>u_{1}(a_{2}); \text{ } u_{2}(a_{1})>u_{2}(a_{2})>u_{2}(c)$};
\node (a) at (1,-4) {$t$};
\node (a) at (2.5,-4) {$f$};
\node (a) at (1.75,-3.5) {$t$};
\node (a) at (0,-6) {$f$};
\node (b) at (1,-6) {$a_2\;a_2\;a_3$};
\draw (0.25,-6.5) rectangle (1.75,-5.5);
\node (c) at (2.5,-6) {$a_3\;a_3\;a_3$};
\draw (1.75,-6.5) rectangle (3.25,-5.5);
\node (a) at (0,-5) {$t$};
\node (b) at (1,-5) {$a_1\;a_2\;a_3$};
\draw (0.25,-5.5) rectangle (1.75,-4.5);
\node (c) at (2.5,-5) {$a_1\;a_3\;a_3$};
\draw (1.75,-5.5) rectangle (3.25,-4.5);
\node (a) at (6,-4) {$t$};
\node (a) at (7.5,-4) {$f$};
\node (a) at (6.75,-3.5) {$f$};
\node (a) at (5,-6) {$f$};
\node (b) at (6,-6) {$a_2\;a_2\;a_2$};
\draw (5.25,-6.5) rectangle (6.75,-5.5);
\node (c) at (7.5,-6) {$c\;\;\;c\;\;\;c$};
\draw (6.75,-6.5) rectangle (8.25,-5.5);
\node (a) at (5,-5) {$t$};
\node (b) at (6,-5) {$a_1\;a_2\;a_1$};
\draw (5.25,-5.5) rectangle (6.75,-4.5);
\node (c) at (7.5,-5) {$a_1\;a_1\;a_1$};
\draw (6.75,-5.5) rectangle (8.25,-4.5);
\path[->] (7.75,-4.75) edge [out=135, in=45] (3,-4.75);
\path[->] (6,-4.75) edge [out=45, in=135] (7.5,-4.75);
\path[->] (6.5,-5.75) edge (6.5,-5.25);
\path[->] (2,-5.25) edge (2,-5.75);
\path[->] (2,-6.25) edge [out=-135, in=-45] (0.5,-6.25);
\path[->] (1.5,-6.25) edge [out=-45, in=-135] (6.5,-6.25);
\node (a) at (4,-2.75) {$g_3$};
\node (a) at (-1,-4.5) {$1$};
\node (b) at (-0.5,-3.5) {$2$};
\node (b) at (0.5,-2.75) {$3$};
\path[->] (-1.25,-3.75) edge (-1.25,-4.5);
\path[->] (-1.25,-3.75) edge (-0.5,-3.75);
\path[->] (-1.25,-3.75) edge [out = 75, in = 180] (0.5,-3);
\node (a) at (-0.75,-8) {$g_3^{v_1}$};
\node (a) at (4.25,-8) {$g_3^{v_2}$};
\node (a) at (9.25,-8) {$g_3^{v_3}$};

\node (a) at (-2.5,-9) {$t$};
\node (a) at (-2,-8.5) {$t$};
\node (a) at (-1.5,-9) {$f$};
\node (a) at (-3,-11) {$f$};
\node (b) at (-2.5,-11) {\text{ }$a_2$};
\draw (-2.75,-11.5) rectangle (-2,-10.5);
\node (c) at (-1.5,-11) {$a_3$\text{ }};
\draw (-2,-11.5) rectangle (-1.25,-10.5);
\node (a) at (-3,-10) {$t$};
\node (b) at (-2.5,-10) {\text{ }$a_1$};
\draw (-2.75,-10.5) rectangle (-2,-9.5);
\node (c) at (-1.5,-10) {$a_1$\text{ }};
\draw (-2,-11.5) rectangle (-1.25,-9.5);
\path[->] (-1.8,-11) edge (-2.2,-11);
\path[->] (-1.5,-10.2) edge (-1.5,-10.8);
\path[->] (-2.5,-10.8) edge (-2.5,-10.2);
\node (a) at (-0.5,-9) {$t$};
\node (a) at (0,-8.5) {$f$};
\node (a) at (0.5,-9) {$f$};
\node (a) at (-1,-11) {$f$};
\node (b) at (-0.5,-11) {\text{ }$a_2$};
\draw (-0.75,-11.5) rectangle (0,-10.5);
\node (c) at (0.5,-11) {$c$\text{ }};
\draw (0,-11.5) rectangle (0.75,-10.5);
\node (a) at (-1,-10) {$t$};
\node (b) at (-0.5,-10) {\text{ }$a_1$};
\draw (-0.75,-10.5) rectangle (0,-9.5);
\node (c) at (0.5,-10) {$a_1$\text{ }};
\draw (0,-11.5) rectangle (0.75,-9.5);
\path[->] (-0.2,-11) edge (0.3,-11);
\path[->] (-0.5,-10.8) edge (-0.5,-10.2);
\path[->] (0.5,-10.2) edge (0.5,-10.8);
\path[->] (-1.5,-11.2) edge [out=-45, in=-135] (0.5,-11.2);


\node (a) at (2.5,-9) {$t$};
\node (a) at (3,-8.5) {$t$};
\node (a) at (3.5,-9) {$f$};
\node (a) at (2,-11) {$f$};
\node (b) at (2.5,-11) {\text{ }$a_2$};
\draw (2.25,-11.5) rectangle (3,-10.5);
\node (c) at (3.5,-11) {$a_3$\text{ }};
\draw (3,-11.5) rectangle (3.75,-10.5);
\node (a) at (2,-10) {$t$};
\node (b) at (2.5,-10) {\text{ }$a_2$};
\draw (2.25,-10.5) rectangle (3,-9.5);
\node (c) at (3.5,-10) {$a_3$\text{ }};
\draw (3,-11.5) rectangle (3.75,-9.5);
\node (a) at (4.5,-9) {$t$};
\node (a) at (5,-8.5) {$f$};
\node (a) at (5.5,-9) {$f$};
\node (a) at (4,-11) {$f$};
\node (b) at (4.5,-11) {\text{ }$a_2$};
\draw (4.25,-11.5) rectangle (5,-10.5);
\node (c) at (5.5,-11) {$c$\text{ }};
\draw (5,-11.5) rectangle (5.75,-10.5);
\node (a) at (4,-10) {$t$};
\node (b) at (4.5,-10) {\text{ }$a_2$};
\draw (4.25,-10.5) rectangle (5,-9.5);
\node (c) at (5.5,-10) {$a_1$\text{ }};
\draw (5,-11.5) rectangle (5.75,-9.5);
\path[->] (4.8,-10) edge (5.2,-10);
\path[->] (5.5,-10.2) edge (5.5,-10.8);
\path[->] (3.5,-11.2) edge [out=-45, in=-135] (5.5,-11.2);
\path[->] (5.5,-9.8) edge [out=135, in=45] (3.5,-9.8);
\path[->] (3.2,-10) edge (2.8,-10);
\path[->] (3.2,-11) edge (2.8,-11);
\path[->] (4.8,-11) edge (5.2,-11);


\node (a) at (7.5,-9) {$t$};
\node (a) at (8,-8.5) {$t$};
\node (a) at (8.5,-9) {$f$};
\node (a) at (7,-11) {$f$};
\node (b) at (7.5,-11) {\text{ }$a_3$};
\draw (7.25,-11.5) rectangle (8,-10.5);
\node (c) at (8.5,-11) {$a_3$\text{ }};
\draw (8,-11.5) rectangle (8.75,-10.5);
\node (a) at (7,-10) {$t$};
\node (b) at (7.5,-10) {\text{ }$a_3$};
\draw (7.25,-10.5) rectangle (8,-9.5);
\node (c) at (8.5,-10) {$a_3$\text{ }};
\draw (8,-11.5) rectangle (8.75,-9.5);
\node (a) at (9.5,-9) {$t$};
\node (a) at (10,-8.5) {$f$};
\node (a) at (10.5,-9) {$f$};
\node (a) at (9,-11) {$f$};
\node (b) at (9.5,-11) {\text{ }$a_2$};
\draw (9.25,-11.5) rectangle (10,-10.5);
\node (c) at (10.5,-11) {$c$\text{ }};
\draw (10,-11.5) rectangle (10.75,-10.5);
\node (a) at (9,-10) {$t$};
\node (b) at (9.5,-10) {\text{ }$a_1$};
\draw (9.25,-10.5) rectangle (10,-9.5);
\node (c) at (10.5,-10) {$a_1$\text{ }};
\draw (10,-11.5) rectangle (10.75,-9.5);
\path[->] (9.5,-10.8) edge (9.5,-10.2);
\path[->] (9.8,-11) edge (10.2,-11);
\path[->] (10.5,-10.2) edge (10.5,-10.8);
\path[->] (10.5,-9.8) edge [out=135, in=45] (8.5,-9.8);
\path[->] (8.5,-11.2) edge [out=-45, in=-135] (10.5,-11.2);



\node (a) at (0.25,-13) {$U_3 \text{:}$};
\node (a) at (4,-12.5) {$u_{1}(a_{2})>u_{1}(a_{1})>u_{1}(a_{3})>u_{1}(c)$};
\node (a) at (4,-13) {$u_{2}(a_{3})>u_{2}(a_{2})>u_{2}(a_{1})>u_{2}(c)$};
\node (a) at (4,-13.5) {$u_{3}(a_{1})>u_{3}(a_{3})>u_{3}(a_{2})>u_{3}(c)$};
\end{tikzpicture}
  \caption{Normal forms $g_2$ and $g_3$ 
of the non-initialized game structures  $\cG_2$ and $\cG_3$ 
from Figures \ref{Fig1} and \ref{Fig2}, respectively. Also normal forms  $g^{v_1}_2$ and $g^{v_2}_2$ of $\cG_2$ initialized in $v_1$ and  $v_2$, 
and normal forms $g^{v_1}_3$, $g^{v_2}_3$, $g^{v_3}_3$ of $\cG_3$ 
initialized in $v_1$, $v_2$, $v_3$, respectively.}
\label{Fig3}
\end{figure}


\subsection{Mixed and stationary mixed strategies}
\label{ss02}
A {\em mixed strategy}  $x^i$  of a player  $i \in I$
is defined as a probability distribution over the set
of his pure strategies.
Thus, the dimension of this set is  $k_i - 1$.

\smallskip

A {\em stationary mixed} strategy  $y^i$  
of a player  $i \in I$
is defined as a set of probability distributions 
for all $v \in V_i$:
each one over all moves $(v,v') \in E$  from  $v$.
The moves are chosen randomly,  
in accordance with these probability distributions,  
and {\em independently} for all $v \in V_i$.

The dimension of the set of stationary mixed strategies 
of a player $i$  is equal to    
$$k'_i = \sum_{j = 1}^{t_i} (\ell_j^i -1) = 
\sum_{j = 1}^{t_i} \ell_j^i - t_i$$
Obviously,  $k'_i \leq k_i$   and the equality holds 
if and only if  $|V_i| = 1$.
Thus, by definition, the set of stationary mixed strategies
is a subset of the set of mixed strategies.
For the play-once games, and only in this case, 
the above two sets coincide.

WLOG we assume that there are no forced positions, that is, 
$\ell^i_j > 1$ for all $i \in I$ and $j = 1, \ldots, t_i$. 

\begin{remark}
Stationary mixed strategies are closely related to the so-called
{\em behaviour strategies} introduced in  \cite{Kuh50}
for games with imperfect information; see also \cite{Kuh53,Aum64}.
\end{remark}

\subsection{Markov and a-priori realizations; 
expected payoffs}
\label{ss03}
A non-initialized game structure $\cG = (G,D,q)$ 
defines a probability distribution 
$q(v)$ for each  $v \in V_R$. 
Furthermore, each stationary mixed strategy 
$y^i$ of a player  $i \in I$ is 
a set of probability distributions 
$p(v)$  for each  $v \in V_i$. 
Thus, given a situation $y = (y^i \mid i \in I)$ 
in stationary mixed strategies, 
one obtains a probability distribution over 
the set of moves $(v,v')$ for all positions 
$v \in V_1 \cup \ldots V_n \cup V_R = V \setminus V_T$.

These distributions naturally define a Markov chain on $G$.
For any initial position $v_0 \in V \setminus V_T$
we can efficiently compute the unique limiting distribution
$P_M = \{P_M(c,v_0), P_M(a,v_0) \mid a \in V_T\}$
over the set of outcomes
$A = V_T \cup \{c\} = \{a_1, \ldots, a_m, c\}$; 
see, for example, \cite{KS60}.

The limiting distribution $P_M$  is defined 
as a function of probabilities from the distributions
$\{q(v) | v \in V_R\}$  and  $\{p(v) | v \in V_1 \cup \ldots \cup V_n\}$.
It is important to note that this 
function may be discontinuous
already in the deterministic play-once case,
$V_R = \emptyset$  and  $|V_i| = 1$  for all  $i \in I$.

\medskip

Consider, for example, 
game structures  $\cG_2$  and  $\cG_3$.  
For a positions  $v_i$, denote by  $p_i$  
the probability to stay on the cycle; 
then  $1 - p_i$  is the probability  to terminate in $a_i$; 
here $i = 1,2$  for $\cG_2$  and  $i = 1,2,3$  for  $\cG_3$.   
\smallskip 

For $\cG_2$, if $p_1 = p_2 = 1$, the play will cycle with probability 1 
resulting in the limiting distribution $(0, 0, 1)$ on $(a_1, a_2, c)$. 
Otherwise, if  $p_1 < 1$  or  $p_2 < 1$, 
for the initial positions  $v_1$  and  $v_2$, 
we obtain the following limiting distributions, respectively: 
\begin{equation}
\begin{aligned}
\left(\frac{1-p_1}{1-p_1p_2}, \frac{p_1(1-p_2)}{1-p_1p_2}, 0\right),\\
\left(\frac{p_2(1-p_1)}{1-p_1p_2},\frac{1-p_2}{1-p_1p_2}, 0\right).\label{lim_d2}
\end{aligned}
\end{equation}

For $\cG_3$, if $p_1 = p_2 = p_3 = 1$, 
the play will cycle with probability 1 
resulting in the limiting distribution 
$(0, 0, 0, 1)$ on $(a_1, a_2, a_3, c)$.
Otherwise, if  $p_i < 1$  for some  $i \in \{1,2,3\}$, 
for the initial positions $v_1$, $v_2$, and $v_3$, 
we obtain the following limiting distributions, respectively:
\begin{equation}
\begin{aligned}
\left(\frac{1-p_1}{1-p_1p_2p_3}, \frac{p_1(1-p_2)}{1-p_1p_2p_3}, \frac{p_1p_2(1-p_3)}{1-p_1p_2p_3},0\right),\\
\left(\frac{p_2p_3(1-p_1)}{1-p_1p_2p_3},\frac{1-p_2}{1-p_1p_2p_3}, \frac{p_2(1-p_3)}{1-p_1p_2p_3},0\right),\\
\left(\frac{p_3(1-p_1)}{1-p_1p_2p_3}, \frac{p_1p_3(1-p_2)}{1-p_1p_2p_3}, \frac{1-p_3}{1-p_1p_2p_3} ,0\right).\label{lim_d3}
\end{aligned}
\end{equation}

It is important to note that for $\cG_n$,  
the limiting probability  
$P_M(c,v_0)$ of the cycle, as a function of  $(p_1, \ldots, p_n)$, 
has a discontinuity at point $(1, \ldots, 1)$; 
it is $1$ when $p_1 = \ldots = p_n = 1$ and $0$ 
otherwise, for any   $n \ge 2$. 
Cases $n = 2$ and $n=3$ were considered above.

\medskip

Typically, for solving graphical games 
in stationary mixed strategies the Markov realization is applied; 
see for example, \cite{KFSV09,MO70}.  
In \cite{BGY13,BGMOV18} 
the following alternative approach was suggested.
Suppose a play revisits a position $v \in V \setminus  V_T$. 
Then, the new move in  $v$  
must coincide with the previously chosen one. 

In other words, before the play begins, 
in each position  $v \in V$ a move  $(v,v') \in E$ is chosen 
according to $q(v)$ for $v \in V_R$ and 
to  $p(v)$  for  $v \in V_1 \cup \ldots \cup V_n$.
After this, the play begins in an initial position and
follows these chosen moves until it terminates or cycles.
This rule defines the {\em a-priori realization}, which
differs essentially from the Markov one. 
Under the latter, a move  $(v,v') \in E$ 
from a position  $v \in V \setminus V_T$ 
is also chosen in accordance with 
a distribution  $p(v)$  or  $q(v)$,  
and such  random choice with 
the same distribution is repeated whenever 
the  play returns to $v$,  
but the resulting move itself is not necessarily repeated.
In contrast for the a-priori realization, 
the limiting distribution
$P_{apr} = \{P_{apr}(c,v_0), P_{apr}(a,v_0) \mid a \in V_T \}$ 
over the set of outcomes
$A = V_T \cup \{c\} = \{a_1, \ldots, a_m, c\}$
is well-defined (unique), whenever an initial position
$v_0 \in V \setminus  V_T$  is fixed.
Furthermore, $P_{apr}$  is a continuous function of probabilities
from the distributions  $q(v), v \in V_R$, and
$p(v), v \in V_1 \cup \ldots \cup V_n$.
Indeed, for any play
(terminal one or a lasso) beginning in  $v_0$,
its probability equals the product of
probabilities of all moves involved in this play.
Then, to compute $P_{apr}(c,v_0)$  and  $P_{apr}(v, v_0)$ 
for all  $ v \in V_T$, 
we ``simply" sum up the probabilities of the corresponding plays: all lassos in the former case
and all  plays terminating in  $v$  in the latter case.

\medskip 

Let us note, however, that 
the number of plays may be exponential 
in the size of digraph  $G$.
So, unlike the Markov case, the above algorithm
computing the limiting distribution is not efficient.
Whether a polynomial one exists is an open problem.
We conjecture that it does not.

\begin{remark}
We acknowledge that the Markov realization 
has many advantages with respect to 
(WRT) the a-priori one:
the Markov one has many practical applications, 
the limiting distribution can be efficiently computed, etc.
Yet, we will show that
Nash-solvability in mixed strategies
of an initialed play-once game 
holds under the a-priori realization.   
In contrast, under the Markov one,  
Nash-solvability may fail 
(already in the  play-once case)  
because of the discontinuity 
of the Markov expected payoff; 
see Section 2.1  and also 
\ref{s4} for more details.

After all, when a play revisits a position, 
why should we roll the dice again? 
We already did it and can reuse the result. 
\end{remark}

Let us consider  $\cG_2$. 
Assuming that the initial positions are $v_1$ and $v_2$, 
we obtain the following limiting a-priori distributions 
for the outcomes $(a_1, a_2, c)$, respectively:
\begin{equation}
\begin{aligned}
(1-p_1, p_1(1-p_2), p_1p_2),\\
(p_2(1-p_1), 1-p_2, p_2p_1).\\
\label{lim_d2_ap}
\end{aligned}
\end{equation}
For $\cG_3$, assuming that the initial positions 
are  $v_1$, $v_2$, $v_3$, we obtain the following 
limiting a-priori distributions, 
for the outcomes $(a_1, a_2, a_3, c)$, respectively:
\begin{equation}
\begin{aligned}
(1-p_1, p_1(1-p_2), p_1p_2(1-p_3), p_1p_2p_3),\\
(p_2p_3(1-p_1), 1-p_2, p_2(1-p_3), p_2p_3p_1),\\
(p_3(1-p_1), p_3p_1(1-p_2), 1-p_3, p_3p_1p_2).
\label{lim_d3_ap}
\end{aligned}
\end{equation}

The probability of outcome $c$ is $p_1p_2$ for $\cG_2$ 
and $p_1p_2p_3$  and $\cG_3$, 
and it is strictly positive whenever 
$p_i > 0$, for all $i \in I$. 
Indeed, in contrast to the Markov realization, 
under the a-priori one, the cycle will be repeated infinitely 
whenever it appears once.

\medskip

Given a limiting distribution
$P = \{P(c,v_0), P(a,v_0) \mid a \in V_T \}$ over the set of outcomes $A$, 
which can be the Markov or the a-priori one, 
and a payoff function  
$u : I \times A \rightarrow \mathbb{R}$, 
the expected payoff is defined as the linear combination:
$$u_i(P) = P(a_1, v_0) u_i(a_1) + \ldots + P(a_m, v_0) u_i(a_m)  + P(c, v_0) u_i(c),$$
$i \in I, v_0 \in V \setminus V_T.$

\smallskip

When  $P$ is defined by a situation 
$y = (y^i \mid i \in I)$  
in mixed stationary strategies, 
we will use the notation  
$u = (u_i(y) \mid i \in I)$.

\subsection{Nash equilibria  and  Uniform Nash equilibria 
in pure, mixed, and stationary mixed strategies}
\label{ss05}
\subsubsection*{Normal form games}
Given a set of players  $I =   \{1, \ldots, n\}$,
a finite set  $S^i$  of 
pure strategies of each player $i \in I$,
and a set  of outcomes  $A$,
a {\em game form} is defined as a mapping
$g : S \rightarrow A$  that  assigns an outcome  $a \in A$
to each situation 
$s = (s^1, \ldots, s^n) \in S = S^1 \times \ldots \times S^n$.

Given also a payoff function  
$u : I \times A \rightarrow \R$,
where  $u(i,a) = u_i(a)$  
is interpreted as a profit of player $i \in I$
in case of outcome  $a \in A$,
the pair $(g, u)$  defines a {\em game in normal form}.

A situation  $s \in S$  is called
a {\em Nash equilibrium (NE) in pure strategies} 
in game  $(g,u)$  if
$u_i(g(s)) \geq  u_i(g(s'))$  for any  $i \in I$
and for any situation  $s' \in S$
that may differ from  $s$  only in the components  $i$.
In other words,  $s \in S$ is a  NE
if no player $i \in I$ can improve $s$  for himself
by choosing some strategy, $s'^i$  instead of  $s^i$,
provided other players $(j \in I \setminus \{i\})$
apply their old strategies  $s^j$. 

In this case, $s^i$ is called a {\em best response} 
to the strategies $(s^j \mid j \in I \setminus \{i\})$. 
Thus, situation $s = (s^1, \ldots, s^n)$  is an NE 
if and only if  
the strategy of each player is a best response 
to the strategies of the remaining $n-1$  players. 

\begin{remark}
Frequently, game  $(g, u)$  has no NE in pure strategies.
However, some conditions on game form  $g$  
may guarantee Nash-solvability.
For  $n = 2$  such necessary and sufficient conditions
were obtained in \cite{Gur75,Gur89}; see also \cite{GN22}.
However, these conditions do not work for  $n > 2$ \cite{Gur89}; 
see also  \cite{BG03,BGMOV18}.
\end{remark}

\medskip

A {\em mixed strategy} $x^i \in X^i$  
of a player  $i \in I$
is defined as a probability distribution over $S^i$ 
determined by probabilities  
$p(s^i \mid x^i)$  for all $s^i \in S^i$.
Then, each situation
$x = (x^1, \ldots, x^n) \in X^1 \times \ldots \times X^n = X$ in mixed strategies uniquely determines a probability 
distribution over $S$ given by probabilities
$p(s \mid x) = \prod_{i \in I}p(s^i \mid x^i)$  for all
$s = (s^1, \ldots, s^n) \in S$, and 
also expected payoff    
$u : I \times X \rightarrow \R$, where 
$u(i,x) = u_i(x) = \sum_{s \in S} p(s \mid x) u_i(g(s))$  
is the expected payoff of player $i$ in situation  $x$.

\medskip

A situation  $x \in X$  is called
a {\em NE in mixed strategies} 
in the normal form game $(g,u)$ 
if  $u_i(x) \geq  u_i(x')$  for every  $i \in I$
and each situation  $x' \in X$
that may differ from  $x$  only in the component  $i$.

In other words, situation 
$x = (x^i \mid i \in I = \{1, \dots, n\})$ 
is an NE if and only if 
$x^i$  is a best response of player $i$  to 
the strategies $(x^j \mid j \in I \setminus \{i\})$ 
of the remaining $n-1$ players.

Nash \cite{Nas50,Nas51} proved that 
every normal form game  $(g,u)$
has an NE in mixed strategies. 

\subsubsection*{Graphical games}
We apply the above definitions
to the normal form of
an initialized graphical game 
$(\cG, v_0, u)$  
to obtain the following concepts:
\begin{enumerate}[(i)]
\item  NE in pure strategies;
\item NE in stationary mixed strategies under the Markov realization;
\item NE in stationary mixed strategies under the a-priori realization.
\end{enumerate}

Note that the concept of an NE in mixed strategies
can also be defined by the mixed extension
of the normal form, yet, it is not realized
by an position-wise independent randomization.
However, mixed and stationary mixed  
strategies coincide for the play-once games.

\subsubsection*{Uniform Nash equilibria}
Given a non-initialized graphical game $(\cG, u)$,
we define a {\em uniform NE} (UNE) 
as a situation which is an NE  in  $(\cG, v_0, u)$  
for any  $v_0 \in V \setminus V_T$.
This modification is applicable 
in all three above cases: (i), (ii), and (iii).
Two examples of the UNE-free games 
are given in Figures 1 and 2.

\begin{remark}
The name of {\em subgame perfect} 
NE is common in the literature,
but we prefer to call such NE uniform,
because the concept of a subgame itself 
loses its meaning in presence of cycles.
\end{remark}

\section{Main results}
\label{s1}

\subsection{Markov realization}
\label{ss11}

\subsubsection*{Uniform best responses}
Note that a Markov decision process can be 
viewed as a graphical one-person game under the Markov realization.
The main result in this area states that
there exists a uniform best pure strategy,
which can be found as a solution of a linear program \cite{Hov60,MO70}.
(As usual, ``uniform best" means 
``best WRT any initial position 
$v \in V \setminus V_T$".)
For the $n$-person case this result 
can be reformulated as follows.

\begin{proposition}
\label{p1}
Given an $n$-person graphical game 
under the Markov realization,
for any set of mixed stationary strategies 
$(y^j \mid j \in I \setminus \{i\})$
of $n-1$  players there exists 
a uniform best response of player  $i$ 
in pure strategies.
\qed
\end{proposition}

\begin{remark} 
This result in a more general setting, 
namely, for stochastic games 
with countable state and action spaces 
and semi-Markov strategies was obtained in \cite{HKW83}.  
The proof is based 
on the results of \cite{Bla62, Sha53}.
\end{remark} 

\subsubsection*{On Nash equilibria in
initializing extensions of graphical games}
In its turn, the last statement implies 
the following relation
between UNE in stationary mixed strategies
in a non-initialized game structure  $(\cG, u)$  
and NE in its initializing extension  $(\cG', u)$;  
see Section 1.1 for definitions 
and Figures \ref{Fig1} and \ref{Fig2} for examples.

\begin{proposition}
\label{p2}
Given a situation $y$ in stationary mixed  strategies
in a non-initialized game  $(\cG,u)$,

(i) if $y$  is an UNE in  $(\cG, u)$ then
$y$  is an NE in  $(\cG', u)$
WRT every distribution  $q(v_0)$.

(ii) if  $y$  is an NE in $(\cG',u)$
for some strictly positive $q(v_0)$  then
$y$  is an UNE in $(\cG, u)$.
\end{proposition}

This statement appeared in \cite{BEGM12} 
for the case of pure stationary strategies. 
Here we extend it 
to the case of stationary mixed strategies.

\smallskip

\proof
Implication  (i) $\Rightarrow$  (ii) is obvious.
If  $x$  is an NE  $(\cG, u)$
WRT every initial position  $v \in V \setminus V_T$  then
$x$  is an NE in $(\cG', u)$  WRT  $v_0$.
Indeed, $v$  will follow  $v_0$  
with probability  $q(v_0,v)$ and
after this play never returns to  $v_0$.
Hence, all expected payoffs in $(\cG', u)$  
initialized in  $v_0$
equal linear combinations of the corresponding 
payoffs in  $(\cG, u)$  initialized  in $v$
with non-negative coefficients  $q(v_0,v)$.
This operation respects inequalities.

\medskip

Implication  (ii) $\Rightarrow$ (i) 
follows from Proposition \ref{p1}.
Suppose that  $x = (x^i \mid i \in I)$  
is not a UNE in  $(\cG, u)$.
Then, there is a position  
$v^* \in V \setminus V_T$  and a player $i \in I$
who can strictly improve situation $x$  
for himself replacing  $x^i$  by  $\tilde{x}^i$,
provided the game  begins in  $v^*$.
But $i$  has a uniform best response 
in situation $x$.
WLOG we can assume  that it is $\tilde{x}^i$.
Hence, $\tilde{x}^i$  
strictly improves  $x$  for  $i$
when the game  begins in  $v^*$  and
it gets  at least as good result as in  $x$
when the game begins in any position $v \in V_T$,
just because $x^i$  is also a response in  $x$.
Thus, $x$ is not an NE in  
$(\cG', u)$  provided  $q(v_0,v) > 0$.
\qed

\smallskip

We will see that the last condition is essential.

\subsubsection*{On Nash equilibria 
in games $(\cG_2, u)$  and $(\cG_3, u)$}
It was shown in \cite{AGH10} that game   $(\cG_2, u)$
has no NE in pure strategies if and only if $u \in U_2$,
where $U_2$ is defined by the system of inequalities:

\begin{equation}
    u_{1}(c)>u_{1}(a_{1})>u_{1}(a_{2}); \text{ } u_{2}(a_{1})>u_{2}(a_{2})>u_{2}(c).\label{U2}
\end{equation}

In 
\ref{ss41},  
we will extend this result 
to the case of mixed strategies as follows. 

\begin{proposition}
\label{p3}
Game  $(\cG_2, u)$  has no UNE 
in mixed strategies
when $u \in U_2$.
\end{proposition}

Then, by Proposition \ref{p2}, 
the following statement holds.

\begin{proposition}
\label{p3_}
Initializing extension  $(\cG'_2, u)$ 
has no NE in mixed strategies when  $u \in U_2$.
\qed 
\end{proposition}

\bigskip

It was shown in \cite{BG03} that game $(\cG_3, u)$
has no UNE in pure strategies if  $u \in U_3$,
where $U_3$  is defined by the following system of inequalities:
\begin{equation}
\begin{aligned}
u_{1}(a_{2})>u_{1}(a_{1})>u_{1}(a_{3})>u_{1}(c),\\ u_{2}(a_{3})>u_{2}(a_{2})>u_{2}(a_{1})>u_{2}(c),\\
u_{3}(a_{1})>u_{3}(a_{3})>u_{3}(a_{2})>u_{3}(c).\label{U3}
\end{aligned}
\end{equation}

In 
\ref{ss41},  
we strengthen this result as follows. Let us set 
\begin{equation}
    \mu_1 = \frac{u_{1}(a_{2})-u_{1}(a_{1})}{u_{1}(a_{1})-u_{1}(a_{3})},
    \mu_2 = \frac{u_{2}(a_{3})-u_{2}(a_{2})}{u_{2}(a_{2})-u_{2}(a_{1})},
    \mu_3 = \frac{u_{3}(a_{1})-u_{3}(a_{3})}{u_{3}(a_{3})-u_{3}(a_{2})}.\label{mu}
\end{equation}
Obviously, $\mu_i > 0$ for $i = 1, 2, 3$ when $u \in U_3$.
\begin{proposition}
\label{p4}
For $u \in U_3$, game  $(\cG_3, u)$  has no UNE in mixed strategies
when  $\mu_1 \mu_2 \mu_3 \geq 1$. 
Otherwise, if $\mu_1 \mu_2 \mu_3 < 1$, 
game $(\cG_3, u)$ has a unique UNE 
in mixed strategies determined by probabilities
\begin{equation}
\begin{aligned}
p_1 = \frac{\mu_3(1+ \mu_1+\mu_1\mu_2)}{1+\mu_3 + \mu_3\mu_1},\\
p_2 = \frac{\mu_1(1 + \mu_2+\mu_2\mu_3)}{1+\mu_1 + \mu_1\mu_2},\\
p_3 = \frac{\mu_2(1 + \mu_3+\mu_3\mu_1)}{1+\mu_2 + \mu_2\mu_3}.
\label{p_i}
\end{aligned}
\end{equation}
\end{proposition}

This statement appears in \cite{BGY13},
yet, no complete proof was given; 
we will give it in 
\ref{s4}.
Right now let us only note that  
each of the three equations of \eqref{p_i} implies that   
$\mu_1\mu_2\mu_3 \leq 1$ and 
this inequality is strict whenever 
$p_i<1$, $i=1, 2, 3$.

\medskip

By Proposition \ref{p2}, we conclude that
the initializing extension  $(\cG'_3, u)$  of this game
has the same NE in mixed strategies, 
provided  $u \in U_3$ and all probabilities  
$q(v_0,v)$  are strictly  positive. 
The last condition is essential. 
It is not difficult to verify that 
if $u  \in U_3$ and 
$q(v_0) = (q(v_0,v_1),( q(v_0,v_2), 
q(v_0,v_3)) = (1,0,0)$ or $(1/2,0,1/2)$
then  $p = (p_1, p_2, p_3) = (0,1,1)$ is 
a pure strategy NE in game $(\cG'_3, u)$,
while  $(\cG_3, u)$  has no UNE when  $u \in U_3$.

However, if we restrict the players 
to their strictly mixed 
strategies ($p_i > 0$ for $i = 1, 2, 3$), 
then games $(\cG'_3, u)$  for all $q(v_0)$ 
become equivalent to  $(\cG_3, u)$, 
that is, all these games have the same NE.

\begin{proposition}
\label{p5_}
For  $u \in U_3$, 
game  $(\cG'_3, u)$  has a unique 
mixed strategy NE given by \eqref{p_i} 
when  $\mu_1 \mu_2 \mu_3 < 1$; otherwise, if  
$\mu_1 \mu_2 \mu_3 \ge 1$, then game $(\cG'_3, u)$ 
has no NE in strictly mixed strategies.
\end{proposition}

\medskip

Propositions \ref{p3}, \ref{p4}, and \ref{p5_} 
will be proven in   
\ref{ss42}.

\medskip 

Finally, let us recall that 
mixed and stationary mixed strategies coincide
for the play-once game structures
$\cG_2, \cG_3$   and  $\cG'_2, \cG'_3$.

\subsubsection*
{Why does Nash's theorem fail 
in case of a Markov realization?}
Indeed, at the first glance, one may decide that
an NE in mixed strategies must exist in games
$(\cG'_2, u)$  and  $(\cG'_3, u)$
(and more generally, in the initializing extension
of every play-once game)
due to the classical Nash theorem \cite{Nas50,Nas51}.
Yet, it works only in case of the a-priori realization, 
but not for the Markov one.
As we have already mentioned, 
in the latter case the limiting distribution
(and hence the expected payoff as well)
may be a discontinuous function 
of probabilities $p(v)$ and  $q(v)$.

\begin{remark}
Although there are results that extend
Nash-solvability to some discontinuous 
payoff functions 
(see, for  example, \cite{DM86,Ros65})  yet,
these results do not cover the Markov realization.
\end{remark}

\subsection{A-priori realization}
Nash's theorem implies existence  of an NE in every
initialized play-once graphical game
under the a-priori realization.
One obtains such an NE just solving in mixed strategies
the normal form of this game.

Naturally, Nash-solvability in mixed stationary strategies
may fail if the game is not play-once;
see, e.g., the main example in \cite{BGMOV18}.
This is not a surprise, since 
the mixed and stationary mixed strategies 
coincide only for the play-once games, 
otherwise the former set is 
a proper subset of the latter.

 \medskip

The {\em uniform} Nash-solvability 
may fail even in the play-once case. In 
\ref{ss42}, we show it for game structures 
$\cG_2$  and  $\cG_3$.

\begin{proposition}
\label{p6}
Under the a-priori realization, games  $(\cG_2, u)$  and  
$(\cG_3, u)$ have no UNE in mixed strategies
whenever  $u \in U_2$  and  $u \in U_3$, respectively.
\end{proposition}


Moreover, a uniform best strategy may fail to exist 
already for one player, 
that is, for a Markov decision process.   
Consider, for example, game structure $\cG_1$ in Figure \ref{Fig4}.

\begin{figure}[t]
\begin{center}
\begin{tikzpicture}
\node[ell] (e1)at (-10,0) {$v_1$};
\node[ell] (e2)at (-7,0) {$v_0$};
\node[ell] (ea1)at (-10,2) {$a_1$};
\node[ell] (ea2)at (-7,2) {$a_2$};
\node[] (through2) at (-8.5, -0.5) [below]{\footnotesize$1/2$};
\node (node1) at (-11,3) {$\cG_1:$};

\draw [->] (e1) to node[below]{\footnotesize$p$} (e2);
\draw [->] (e1) to node[left]{\footnotesize$1-p$} (ea1);
\draw [->] (e2) to node[right]{\footnotesize$1/2$} (ea2);
\draw[->,rounded corners] (e2.south) |- (through2.south) -| (e1.south);
\end{tikzpicture}
\end{center} \label{Fig0}
\caption{
Game structure  $\cG_1$ 
in which $v_0$ is a position of chance 
and $v_1$ is a position of player 1.}
\label{Fig4}
\end{figure}

The player controls position $v_1$, while $v_0$
is a position of chance with two equal probabilities:
$1/2$ and $1/2$.
Consider any  payoff  $u$  satisfying inequalities
$$\frac{1}{2}\left(u(a_2) + u(c)\right) > u(a_1) > u(c).$$
Then, if game begins in $v_1$
the optimal player's move is to  $v_0$,
while if the initial position is  $v_0$,
it is better to terminate in  $a_1$, avoiding  $c$.
Note that this happens only under the a-priori realization,
while under the Markov realization, move $(v_1,v_0)$
will be the best for both initial position: 
$v_0$ or  $v_1$.

This argument shows why Proposition \ref{p2}
cannot be extended to the a-priori realization:
compare Proposition \ref{p4} and
the opening  claim of this  subsection.

\section{Two main examples are UNE-free}
\label{s3}
Our two main examples are given 
by two play-once non-initialized
2- and 3-person game structures $\cG_2$  and  $\cG_3$ 
given on Figures \ref{Fig1} and \ref{Fig2}, respectively.   
In normal form both are represented in Figure \ref{Fig3}. 
The corresponding games $(\cG_2, u)$  and  $(\cG_3,u)$
have no UNE in pure strategies whenever 
payoffs are ordered in accordance 
with systems of strict inequalities  $U_2$  and  $U_3$
defined by \eqref{U2} and \eqref{U3}, respectively.
Interpretation of both games are given in 
\ref{s7}. 

\subsection{Game $\cG_2(u)$  for  $u \in U_2$}
\label{ss32}
Game $(\cG_2, u)$ has no UNE when $u \in U_2$.

The game is play-once. 
Each player $i \in I = \{1,2\}$  
controls a unique  position  $v_i$
and has two strategies:
either to terminate in  $a_i$  or to follow the cycle: 
$s^i \in \{t, f\}$.

We have to show that none of the four situations
$s = (s^1, s^2)$ is a UNE, that is,
at least one player can improve it
WRT at least one initial position
$v_0 =  v_i, \; i \in \{1,2\}$. 

\smallskip

Consider  $s = (f,f)$:  
all players follow the cycle.
The play results in  $c$  
for any initial position.   
Player 2  can improve his result
choosing  $t$  rather than  $f$  and getting  $a_2$, 
WRT any initial position.

\smallskip

Consider  $s = (f,t)$.
The play results in  $a_2$  
WRT  any initial position.  
Player 1 can improve her result 
WRT  $v_1$ 
choosing  $t$  rather than  $f$  
and terminating in  $a_1$  instead of  $a_2$. 
Yet,  WRT  $v_2$ there is no improvement. 

\smallskip

Consider  $s = (t,t)$.
The play results in  $a_i$
WRT initial position  $v_i$, for $i = 1,2$. 
Player 2 can improve his result 
WRT initial position  $v_2$ 
choosing  $f$  rather than  $t$  
and terminating in  $a_1$  instead of $a_2$. 
Yet,  WRT  $v_1$ the outcome is  $a_1$ 
for both his strategies. 

\smallskip

Consider  $s = (t,f)$.
The play results in  $a_1$
for any initial position. 
Player 1 can improve her result 
WRT any initial position  
choosing  $f$  rather than  $t$  
and getting  $c$  instead of $a_1$.  

\smallskip 

Thus, we obtain  $s = (f,f)$  again. 
All four situations belong to 
an improvement cycle of length 4. 
It is shown in Figure \ref{Fig3}. 
Hence, none of them is a UNE 
in game  $(\cG_2, u)$  when  $u \in U_2$. 

\medskip 

In contrast, both initialized game structures 
$\cG^{v_1}_2$  and  $\cG^{v_2}_2$, 
with initial positions  $v_1$  and  $v_2$, 
as well as the corresponding game forms 
$g^{v_1}_2$  and  $g^{v_2}_2$, are Nash-solvable; 
see Figures \ref{Fig1} and \ref{Fig3}. 

\medskip 

In fact, Nash-solvability holds for any 
deterministic initialized game structure of two players.
This result was derived in \cite{BG03} 
from a general criterion of Nash-solvability 
\cite{Gur75,Gur89}; see also \cite{GK18,GN22}. 
However, this criterion holds only for 
$n = |I| = 2$; see Section \ref{s5} for more details.

\subsection{Game $(\cG_3,u)$ for $u \in U_3$}
\label{ss33}
Game $(\cG_3,u)$ has no UNE when $u \in U_3$ \cite{BG03}. 
For completeness, we provide here a simplified proof.  

The game is play-once. 
Each player $i \in I = \{1,2,3\}$  
controls a unique  position  $v_i$
and has two strategies:
either to terminate in  $a_i$  or to follow the cycle;
$s^i \in \{t, f\}$.
We have to show that none of the eight situations
$s = (s^1, s^2, s^3)$ is a UNE when  $u \in U_3$, 
that is, at least one player can improve  $s$ 
WRT at least one initial position 
$v_0 =  v_i, \; i \in \{1,2,3\}$.

\smallskip

Consider  $s = (f,f,f)$:   
all 3 players follow the cycle.
For any initial position 
the play results in  $c$  and
each player can improve his result
choosing  $t$  rather than  $f$,
WRT any initial position.

\smallskip

Consider  $s = (t,t,t)$: 
all 3 players terminate.
Each one can improve the situation 
choosing  $f$  rather than  $t$.
Then, the next player will terminate,
which is better, according to  $U_3$.
Note, yet, that improvement for player  
$i$ is strict only when  $v_0 =  v_i$.
Otherwise, the  outcome will  not  change. 

\smallskip

The remaining six situations
form an improvement cycle. 

Indeed, in situation  
$(f,f,t)$  player  $1$  is unhappy and
will  switch from  $f$  to  $t$.
Doing so (s)he improves the situation, 
at least when $v_0 = v_1$.
In this case  $a_3$  is replaced by  $a_1$,
which is better to  player $1$, according to $U_3$.
Note that the outcome will remain unchanged
when  $v_0 = v_2$  or  $v_0 = v_3$.

The obtained situation $(t,f,t)$  
can be improved by player $3$ 
by switching from  $t$  to  $f$, 
at least when $v_0 = v_2$ or $v_0 = v_3$.
In both cases $a_3$ is replaced by $a_1$,
which is better to  $3$  according to $U_3$.
However, if $v_0 = v_1$, the outcome will not change.

The obtained situation  
$(t,f,f)$ is a ``shift" of $(f,f,t)$, 
which was already considered.
Repeating the same arguments two more times,
we obtain the improvement cycle of length 6:   
$$(f,f,t), (t,f,t), (t,f,f), (t,f,f), (t,t,f),(f,t,f), (f,t,t)$$. It is shown in Figure \ref{Fig3}.
Thus, none of eight situations of game   
$(\cG_3,u)$  is a UNE  in pure stationary strategies, when  $u \in U_3$; 
see Figures 2 and 3.

\medskip 

In contrast, all three initialized game structures 
$\cG^{v_1}_3, \cG^{v_2}_3, \cG^{v_3}_3$  
with initial positions  $v_1, v_2, v_3$, 
as well as the corresponding game forms 
$g^{v_1}_3, g^{v_2}_3$, and $g^{v_3}_3$  
are Nash-solvable; 
see Figures \ref{Fig2} and \ref{Fig3}.

\medskip 

Let us note, however, that 
initialized deterministic $n$-person games 
without NE in pure stationary strategies exist for $n > 2$.  
First such example for $n=4$  
was obtained in \cite{GO15}. 
Then in \cite{BGMOV18} a much smaller 
$3$-person game was constructed 
that has no NE even in stationary mixed strategies.  
However, these NE-free games are not play-once.
It remains an open question, 
whether a play-once NE-free example exist;  
see Section \ref{s5} for more details.

\subsection{Generalizations and possible applications} 
\label{ss13}
A tedious but routine case analysis 
allows to verify that a UNE, in pure stationary strategies, 
exists in games $(\cG_2, u)$  and  ($\cG_3, u)$ 
whenever  $u \not\in U_2$  and  $u \not\in U_3$, respectively.

\medskip 

For  $n > 3$  there are very many similar UNE-free examples; 
some of them will be given in  
\ref{sA1}. 

\smallskip 

Using these examples one can try 
to solve an important open problem:
Construct an initialized deterministic $n$-person game
that has no NE in pure stationary strategies 
and satisfies the following condition

\smallskip

{\bf (C)} \;\; outcome $c$ is worse 
than each terminal outcome $a_j \in A$
for every player  $i \in I$.

\smallskip

Note that $U_3$  satisfies {\bf (C)} while $U_2$  does not.

For more details see Section \ref{s5} and \cite{BGMOV18}, 
in particular, Remark 3 there.
As we already mentioned, for  $n=2$, by \cite{Gur75,Gur89,BG03},
an NE in pure strategies exists, 
even if condition {\bf (C)} is waved.

\section{Main results and open problems related to Nash-solvability 
in pure stationary strategies}
\label{s5}

\subsection{Uniform Nash-solvability 
in presence of moves of chance}  
There are two important classes of games
that always have a UNE in pure stationary strategies:

\begin{itemize}
\item[]
Graphical $n$-person games on acyclic digraphs. 
In this case a special UNE in pure strategy strategies can be found by
Backward Induction \cite{Gal53,Kuh53}; see also \cite{Gur17,Gur17a}.
\item[]
Two-person zero-sum graphical games.
In this case, the existence of an UNE 
follows from basic results of the stochastic game theory
\cite{Gil57,LL69}.
\end{itemize}

More details can be found in  \cite{BG09}. 
In both cases, the uniform Nash-solvability holds 
not only for the terminal effective payoffs,
considered in this paper, but also 
for a wide family of more general types: 
limiting mean, total, or $k$-total effective payoffs; 
see, for example, \cite{BEGM17,EM79,Gil57,GZ04,GKK88,LL69,Mou76,Sha53,TV87,TV98}. 

\medskip 

Let us mention also that 
2-person deterministic graphical games are Nash-solvable; 
see Subsection \ref{ss53}  below. 
Yet, as we know, such games may be NE-free, 
and hence, Nash-solvability may fail 
for initialized non-deterministic graphical games  
with only one (initial) position of chance. 

\subsection{NE-free graphical games with 
a unique position of chance and a unique directed cycle}
Recall our main two examples  
$(\cG'_2,u), u \in U_2$  and  $(\cG'_3,u), u \in U_3$.  
Both game structure  $\cG_2$  and  $\cG_3$     
are play-once, contain a unique directed cycle 
and a unique position of chance, 
which is the initial position, 
and both are not Nash-solvable:    
have no NE in pure stationary strategies 
when $u \in U_2$ and $u \in U_3$, respectively. 
Furthermore, in both cases, by deleting the initial position, 
we obtain a UNE-free non-initialized graphical game, 
$(\cG_2,u), u \in U_2$  and  $(\cG_3,u), u \in U_3$. 

Condition {\bf (C)} of Subsection \ref{ss13} holds 
for $U_3$   but not for $U_2$.  

However, 2-person UNE-free graphical games satisfying 
{\bf (C)} also exist. 
An example $(\cG_6, u)$  with  $u \in U_6$)   
was constructed first in \cite{BEGM12},  
where  $U_6$  is determined.   
The 2-person deterministic game structure $\cG_6$ 
still have only one directed cycle  $C_6$.  
Yet, this game is not play-once: 
players  1  and 2  alternate in  $C_6$, 
so each of them controls three positions. 

It was proven in the present paper that 
games  $(\cG_2,u)$  and  $(\cG_3,u)$   
remain NE-free even in mixed stationary strategies 
under the Markov realization, 
for all  $u \in U_2$  and for some  $u \in U_3$. 
Also $(\cG_6, u)$  remains UNE-free 
for both Markov and a-priori realizations for all  $u \in U_6$; 
a proof was sketched in \cite{BGY13}. 

\smallskip 

Note, however, that corresponding 
initialized game structures   
$\cG'_2$  and $\cG'_3$  
are Nash-solvable for any payoff    
$u$  under the a-priori realization. 
This  follows from the classic Nash theorem 
\cite{Nas50,Nas51},     
which is applicable in case of the a-priori realization, 
because both $\cG_2$  and $\cG_3$  are play-once. 

\medskip 

Thus, already one (initial) position of chance 
may destroy Nash-solvability, 
even for play-once games and for games satisfying {\bf (C)}. 

\medskip 

So, for the rest of this section, we restrict ourselves 
to the so-called {\em deterministic graphical (DG) games}, 
(without positions of chance) and show 
(or sometimes conjecture) that Nash-solvability of such games, 
in pure stationary strategies,  
can be saved by some additional assumptions. 
By default, we assume that considered DG games are initialized 
unless it is explicitly said otherwise. 

\medskip 

For the beginning, let us note that DG games may be NE-free 
under the above assumptions, 
The first example, with  $n=4$,   
was generated by a computer code  \cite{GO15}. 
Then, a much simpler 3-person DG game was constructed 
in \cite{BGMOV18}, where it was also shown that 
this game has no NE not only in pure 
but also in stationary mixed strategies, 
under both the Markov and a-priori realizations.
Yet, this game is not play-once; 
there is player who controls two (adjacent) positions.  

\subsection{Nash-solvable deterministic graphical $n$-person games} 
\label{ss53}
\subsubsection*{Two-person case, $n=2$.}
Nash-solvability of the 2-person DG games 
was derived in \cite{BG03}  from 
Nash-solvability of the so-called {\em tight} game forms. 
The latter result is old.    
For the zero-sum case it was obtained by 
Edmonds and Fulkerson in 1970 \cite{EF70}, see also \cite{Gur73}. 
Then, it was extended to the general case in \cite{Gur75,Gur89}. 
Recently, a much shorter proof was given in \cite{GN22}. 

Let us underline that condition {\bf (C)} is not required 
for Nash-solvability of the 2-person DG games. 

\smallskip 

Although the concept of tight game forms 
can be naturally extended to the case  $n \geq 2$, yet,  
for  $n > 2$  tightness is no longer related to Nash-solvability: 
it is neither necessary \cite{Gur89}, 
nor sufficient \cite{Gur75,Gur89}; see also \cite{BG03}.   
Several new classes of tight game forms were recently found in 
\cite{Gur17,GK18,GN22a,GN22b}. 

In \cite{Gur17} 
the class of the DG games is extended 
to a larger class of the so-called multi-stage DG games. 
The outcomes of a DG game are formed by all 
terminals of its digraph  $G$  
and one special outcome  $c$  
corresponding to all infinite plays of  $G$. 
In contrast, the outcomes of a multi-stage DG game 
are formed by the strongly connected components of its digraph;  
furthermore, some outcomes may be merged. 
It is shown in \cite{Gur17} that 
multi-stage DG game forms are tight. 
This statement is stronger 
than tightness of the DG game forms shown in \cite{BG03}.

\subsubsection*{Play-once DG games and DG games satisfying {\bf (C)}.}
The play-once  $n$-person DG games satisfying {\bf (C)}  are Nash-solvable. 
This is the main result of \cite{BG03}. 
Moreover, we conjecture that 
each of these two conditions is sufficient for Nash-solvability. 
We have no example of an $n$-person NE-free DG game 
that is either play-once or {\bf (C)} holds. 
 
 
A stronger version of the second conjecture, 
(called ``Catch 22")  was suggested in \cite{Gur21}:  
In every NE-free $n$-person DG game 
there exist at least two players for each of which 
outcome  $c$  is better than at least 2 terminal outcomes. 
In other words, $c$  cannot be either 
the worst or the second worst 
for all players, and not even for all but one. 

\subsubsection*{Symmetric digraphs.}
The digraph  $G$  is called {\em symmetric} if  
$(v',v'')$  is its arc whenever  $(v'',v')$  is 
unless  $v'$  or  $v''$  is a terminal. 
Recently it was shown in \cite{BFGV22} 
that every $n$-person DG game on a symmetric digraph is Nash-solvabile. 
Condition {\bf (C)} is not needed, although it simplifies the proof. 

\medskip 

A wider class 
of the so-called $n$-person {\em shortest path games} 
was also studied in \cite{BFGV22}. 
A local cost $\ell(i,e) = -r(i, e)$  is defined 
for each player  $i \in [n] = \{1, \dots , n\}$  and 
move  $e$  of such game. 
Condition  {\bf (C)}  holds if all  $\ell(i,e) > 0$.
In this case, given a play  $P$, 
the effective cost of  $P$  for  $i$  is 
the sum of the corresponding local costs, 
$L(i, P) = \sum_{e \in P} \ell(i, P)$, 
if  $P$  is a terminal play and 
$L(i, P) = + \infty$  if  $P$  is an infinite play 
(a lasso, whenever all players apply their pure stationary strategies). 

Nash-solvability of the $n$-person shortest path games 
satisfying  {\bf (C)} on symmetric digraphs 
was proven  in \cite{BFGV22}, 
where it was also conjectured that the last condition 
(symmetry) can be waved if  $n=2$. 
This is the so-called {\em bi-shortest path conjecture} \cite{Gur21a}.
However, an NE-free shortest path game exists  
if  $n = 3$  and the digraph is not symmetric \cite{GO14}.   

\medskip

It was also shown in \cite{BFGV22}  that 
a (non-initialized) DG game has a UNE whenever 
(i) its digraph is symmetric, 
(ii) $n=2$, and (iii)  {\bf (C)}  holds.  .   
Conversely, a UNE may fail to exist if 
at least one of the above three conditions fails. 

\medskip 

Somewhat related results were obtained in \cite{BGMS11}. 
For DG games we assume that all lassos form a unique outcome $c$.
The case when all cycles and terminals form pairwise distinct outcomes
was considered in  \cite{BGMS11}, 
where a criterion of Nash-solvability was obtained
for the 2-person such games on symmetric digraphs.  

\section{Graphical games and stochastic games with perfect information} 
\label{s6}
Here we will show that graphical games can be viewed as 
a special subclass of the stochastic games with perfect information 
and, thus, Nash-solvability of both can be studied,  
for the $n$-person case, simultaneously. 

\subsection{On Nash-solvability of mean payoff games} 
\label{ss61} 
Two-person zero-sum stochastic games 
were introduced in 1953 by Shapley \cite{Sha53}. 
In 1957 Gillette \cite{Gil57} considered 
the subclass of stochastic games 
with zero stop probability, 
introduced limiting mean effective rewards for this case,  
and proved the existence of a UNE 
in mixed stationary strategies.
(For the 2-person zero-sum case an NE is just a saddle point.)  
The proof was far from simple; 
Gillete's approach was based 
on the Hardy-Littlewood Tauberian Theorem and   
all conditions of the latter were accurately verified
(and thus the proof finalized) only in 1969 by 
Ligette and Lippman \cite{LL69}.

Also, Gillette outlined 
the subclass of games with perfect information 
and showed that they can be solved in 
uniform optimal and {\em pure} stationary strategies. 

These games remain of interest even in absence 
of moves of chance, 
when two players control all non-terminal position. 
(Each one is controlled by one player.) 
Such games are called {\em deterministic}; 
2-person zero-sum deterministic-stochastic games with 
zero stop probability, perfect information, 
and the limiting average rewards  
are known as the {\em mean payoff games}. 
They were intensively studied since 1970s \cite{Mou76,EM79,GKK88} 
mostly because of the algorithmic complexity 
of their solution \cite{GKK88}.  
No polynomial algorithm 
for the mean payoff games is still known. 
Recently, a quasi-polynomial one 
was obtained for the so-called {\em parity} games, 
which form a special subclass of the mean payoff games \cite{CJKLS17}. 
However, in the present paper we study Nash-solvability  
rather than polynomial solvability. 

All above definitions can be naturally extended 
from the 2-person zero-sum case to the $n$-person one.  
Thus, we can talk about $n$-person 
mean payoff or stochastic games, 
with or without positions of chance. 

\medskip 

Already 2-person 
(but not zero-sum) mean payoff games 
may have no NE in pure stationary strategies. 
The first NE-free example was given in 1988 
\cite{Gur88}; see also \cite{GKK88}. 
It is constructed on 
the complete bipartite  $3 \times 3$  digraph; 
each player controls 3 positions, 
that is, one part of it, and 
the local rewards are symmetric, that is, the same 
for the moves from  $u$  to  $w$  and  from  $w$ to $u$.  
This game can be interpreted as an 
{\em ergodic extension}  of the corresponding 
$3 \times 3$ bimatrix game \cite{Mou76}. 

The normal form of this game is of size 
$3^3 \times 3^3 = 27 \times 27$  and 
it is an open question whether it has an NE 
in stationary mixed  strategies. 

In \cite{Gur90} 
it was shown that this example is, in a way, minimal: 
Every 2-person mean payoff game 
on a bipartite  $2 \times k$  digraph 
has an NE in pure stationary strategies. 

In \cite{BEGM17}  it was shown that 
the above $3 \times 3$ example 
disproves Nash-solvability 
not only of the mean payoff games, 
but also of a much larger family 
of the so-called $k$-total payoff games 
for any integer nonnegative  $k$. 
Case  $k=0$  is associated with the mean payoffs, 
while $k=1$  is assigned to the so-called total payoffs 
introduced in \cite{TV87,TV98}. 

\subsection{Graphical games can be viewed as 
transition-free mean payoff games} 
\label{s62}
Recall that, by definition, 
all infinite plays of a graphical game 
(and in particular, all lassos in its digraph)   
are equivalent, that is, form a single outcome.  
In contrast, mean payoffs depend on 
the directed cycle of the lasso 
that appears in the game after 
all  $n$  players have chosen their 
pure stationary strategies.  

Graphical games can be viewed as a special 
subfamily of mean payoff games 
(with or without positions of chance). 

Given an $n$-person graphical game $(\cG, u)$  
on a digraph  $G$,  
let us add a loop  
$\ell_v$  to each terminal position  $v \in V_T$  in  $G$ 
and for each player  $i \in [n] = \{1, \dots , n\}$ 
set the local reward  $r(i, \ell_v)$  on  $\ell_v$  
equal to the terminal payoff of  $i$  in  $v$. 
Furthermore, set the local reward  $r(i,e) = 0$  
for any other edge of digraph  $G$ 
and each player $i \in [n]$. 

By this construction, in the obtained game 
all its infinite plays 
(more precisely, all plays that 
do not come to a terminal loop)  
are equivalent, since the corresponding effective payoff 
is  0  for each  player,  while on ``finite" plays  
(that end in terminal loops)  
$n$  players may have arbitrary effective payoffs. 

Obviously, condition  {\bf (C)}  holds 
if and only if  $r(i, \ell_v) < 0$  
for each terminal  $v \in V_T$; 
in other words, if the cost of every terminal is 
positive for each player. 

We can naturally call the obtained 
mean payoff games {\em transition-free}, 
because players do not pay for the moves of the play, 
they pay (or are payed) only in the terminals. 
Obviously, the obtained transition-free mean payoff games 
are equivalent with the original graphical game. 

\medskip 

Thus, two main UNE-free examples  
$(\cG_2, u)$  with  $u \in U_2$  and 
$(\cG_3, u)$  with  $u \in U_3$
of the present paper, as well as 
the 2-person UNE-free example  
$(\cG_6, u)$  with  $u \in U_6$   
from \cite{BEGM12}, provide non-initialized 
UNE-free and transition-free mean payoff games. 
Furthermore, the initialized 
3- and  4-person NE-free examples 
from \cite{BGMOV18}  and  \cite{GO15}, respectively,  
provide initialized NE-free and transition-free mean payoff games.

Thus, results of the present paper can be viewed 
within the framework of further studies of  
Nash-solvability in pure stationary strategies 
of stochastic games with perfect information.  

\medskip 

Finally, let us note that Markov realization 
corresponds exactly to solving 
stochastic games in stationary mixed strategies, 
which is standard, 
while the a-priori realization is a different approach, 
which was not applied to stochastic games yet.  

\subsection{Nash-solvability in pure history-dependent strategies} 
In 1997 Thuijsman and Raghavan \cite{TR97} 
proved Nash-solvability in pure history dependent strategies 
for the mean payoff stochastic games with perfect information. 
As we already mentioned, this class of games contains 
graphical games considered in the present paper. 
For this reason, we are trying to solve them in stationary strategies. 

Let us remark that the result of \cite{TR97}  
holds, in fact, for many other classes of effective payoffs, 
in particular, 
for total \cite{TV87,TV98,GO14} and  $k$-total \cite{BEGM17} ones.

\section*{Acknowledgements}
The authors are thankful to the anonymous reviewer 
for many helpful remarks and suggestions. 
The paper was prepared within the framework
of the HSE University Basic Research Program.

\bigskip

\bigskip

{\bf APPENDIX} 

\appendix

\section{A large family of 
$n$-person deterministic graphical games 
without UNE in pure stationary strategies}
\label{sA1}
Consider the  following
$n$-person play-once non-initialized game structure $\cG_n$.
Given a digraph  $G_n = (V,E)$, where

\smallskip

$V = \{v_1, \ldots, v_n; \; a_1, \ldots, a_n\},$ \

\smallskip

$E = \{(v_1,v_2), (v_2,v_3), \ldots, (v_{n-1},v_n), (v_n,v_1); \;
(v_1,a_1), \ldots, (v_n,a_n)\}$.

\smallskip

For any  $n > 2$  set $I = \{1, \ldots, n\}$  and let
each player $i \in I$ make a move in  $v_i$.
A set of payoffs $U_n$ is defined by the following properties:

\begin{enumerate}[(a)]
\item
For each player  $i \in I$  their own terminal  $a_i$
is worse for them than each of the next  $\lfloor n/2 \rfloor$ terminals, in cyclical order. 
(Among themselves these terminals may be ordered arbitrarily and this order may depend on $i$.)
\item
Among the first $\lfloor (n-1)/2 \rfloor$ of them there is at least one, $a_j$, that is worse than $a_i$ for player $j$.
\item
Finally, condition {\bf (C)} of Section \ref{ss13} holds.
\end{enumerate}

Let us note that for $n=3$  conditions (a) - (c) uniquely define the family of payoffs $U_3$, while for $n=2$ they do not define $U_2$.

\begin{proposition}
\label{p5}
The non-initialized play-once $n$-person game
$(\cG_n, u)$  is
UNE-free whenever  $u \in U_n$.
\end{proposition}

\proof
Each player  $i \in I$ controls a unique position  $v_i$
and, thus, has only two pure strategies:
to terminate at $a_i$  or to follow the cycle:
$s^i \in \{t, f\}$.
We have to prove that any situation
$s = \{s^i \mid i \in I\}$  is not a UNE.
Consider three cases.

\smallskip

Case 0. No player terminates, that is, all choose  $f$.
Then, the play results in the cycle and,
by condition {\bf (C)}, each player can improve
choosing  $t$  rather than  $f$. This holds for any initial position,
$v_0 = v_i,  i \in I$.

\smallskip

Case 1. One player  $i \in I$  terminates,
while all others choose  $f$.
Then, by condition  (b), there exists a player
$j \in \{i+1, \ldots, i+\lfloor \frac{n-1}{2} \rfloor\}$  
who can improve her result
by choosing  $t$  instead of  $f$, at least when $v_0 = v_j$.

\smallskip

Case 2. At least two players terminate.
Obviously, there exist two of them
$i, j \in I$  such that distance from
$v_i$  to $v_j$  along the cycle is at most
$\lfloor \frac{n}{2} \rfloor$.
Then, by  (a), player  $i$  can improve her result 
by switching from  $f$  to  $t$, at least for  $v_0 =  v_i$.
\qed

\bigskip 

\section{Markov and a-priori realizations 
for two main examples; proofs of Propositions 3, 5-7}
\label{s4}
Here we study the uniform Nash-solvability of these games and Nash-solvability of their initializing extensions in the mixed strategies under the Markov and a-priori realizations and prove Propositions \ref{p3}-\ref{p5_}.

\subsection{Games  $(\cG_2, u)$  with  $u \in U_2$  and
$(\cG_3, u)$   with  $u \in U_3$
might have UNE only in strictly mixed strategies,
under both the Markov or a-priori realizations}
\label{ss40} 
As we already know,
games  $(\cG_2, u)$   and   $(\cG_3, u)$
have  no  UNE  in pure strategies  when
$u \in U_2$  and   $u \in U_3$, respectively.
We will strengthen this claim as follows:

\begin{lemma}
\label{l1}
For both the Markov or a-priori realizations,
games  $(\cG_2, u)$  with  $u \in U_2$   and
$(\cG_3, u)$  with  $u \in U_3$
may have  UNE  only in strictly mixed strategies;
in other words,  only when  $0 < p_i < 1$   for $i \in \{1,2,3\}$.
\end{lemma}

\proof
Let  $p = (p_1, p_2)$  be a UNE in  $(\cG_2, u)$  with  $u \in U_2$.
We will show that  if  $p_i$  equals  $0$  or  $1$   then
the same property holds for  $p_{3-i}$, where $i \in \{1,2\}$.
In fact, this was already shown in Section \ref{ss03}:
for every pure strategy of player  $i$
there exists a unique uniform best response of the opponent,
and this response is realized by a pure strategy, while
every strictly mixed response, $0 < p_{3-i} < 1$,
can be improved WRT at least one initial position.

\smallskip

Let  $p = (p_1, p_2, p_3)$  be a UNE in
$(\cG_3, u)$  with  $u \in U_3$.
If  $p_i$  equals  $0$  or  $1$
for a player  $i \in I = \{1,2,3\}$  then
the same property holds for  the two remaining players.
In fact,  we can just repeat the arguments of Section \ref{ss03}.
Since this case is cyclically symmetric
(unlike the previous one)  WLOG  we can set  $i = 3$.

Suppose  $p_3 = 0$, that is, player  $3$  terminates in  $a_3$.
Then player $2$  has a unique uniform best response:
to follow the cycle with move  $(v_2,v_3)$.
Then, player $1$  also has a unique uniform best response:
to terminate with move  $(v_1,a_1)$.

Suppose  $p_3 = 1$, that is, player  $3$  follows the cycle
by move  $(v_3,v_1) $.
Then player $2$  has a unique uniform best response:
to terminate by move  $(v_2, a_2)$.
Then, player $1$  also has a unique uniform best response:
to follow the cycle with move  $(v_1, v_2)$.

In each of the four above cases
the best response is unique and it is realized
by a pure strategy.

\medskip

It is important to note that
\begin{itemize}
\item[]
By definition of an NE  $(p_i \mid i \in I)$,
the strategy  $p_i$  of each player $i \in I$
is a best response (not necessarily unique)
to the set of strategies of the remaining players $I \setminus \{i\}$.
\item[]
All above claims  hold for both the Markov and a-priori realizations.
Although in the latter case a uniform best response
may fail to exist, in general, but for  games,
$(\cG_2, u)$  with   $u \in U_2$,  and
$(\cG_3, u)$  with   $u \in U_3$,
it exists in all considered cases.
\end{itemize}

\begin{remark}
As we know, no UNE in pure strategies exists
for both games under both realizations.
Yet, a UNE in strictly mixed strategies might exist.
This question will be studied in the next two Sections.
\end{remark}

In what follows we denote by $J$ 
the set of indices of non-terminal positions 
and by $F_{ji}$ the expected payoff of player $i$, 
provided the play starts at $v_j$. 
Observe that $F_{ji}$ are continuously differentiable functions 
of $p_i$ when $0< p_i< 1$, $i \in I$. 
Thus, if $p = (p_i\mid i\in I)$ 
is a uniform NE in strictly mixed strategies 
($0< p_i< 1$, $i\ in I$) under 
either the Markov or a-priori realization, then
\begin{equation}
    \frac{\partial F_{ji}}{\partial p_i} = 0, \text{ for all $i \in I, j \in J$}. 
\label{part_deriv}
\end{equation}

\subsection{Markov realization}
\label{ss41}
\begin{proof}[Proof of Proposition \ref{p3}]
Let $p =(p_1, p_2)$ be a uniform NE in game $(\cG_2,u)$ under the Markov realization. If $p_1=p_2=1$, the probability of cycling is 1. Otherwise, the limiting distributions for initial positions $v_1$ or $v_2$ are given by \eqref{lim_d2}, and hence, the expected payoffs are
\begin{equation}
\begin{aligned}
    F_{11} = &\frac{(1-p_1)u_{1}(a_{1})+p_1(1-p_2)u_{1}(a_{2})}{1-p_1p_2},\\
    F_{12} = &\frac{(1-p_1)u_{2}(a_{1})+p_1(1-p_2)u_{2}(a_{2})}{1-p_1p_2},\\
    F_{21} = &\frac{(1-p_2)u_{1}(a_{2})+p_2(1-p_1)u_{1}(a_{1})}{1-p_1p_2},\\
    F_{22} = &\frac{(1-p_2)u_{2}(a_{2})+p_2(1-p_1)u_{2}(a_{1})}{1-p_1p_2}.  \nonumber
\end{aligned}
\end{equation}
In this case, relations \eqref{part_deriv} have the following form:
\begin{equation}
\begin{aligned}
\frac{(u_{1}(a_{1}) - u_{1}(a_{2}))(1-p_2)}{(1-p_1p_2)^2}=0,\\
\frac{-p_1(u_{2}(a_{1}) - u_{2}(a_{2}))(1-p_1)}{(1-p_1p_2)^2}=0.
\nonumber
\end{aligned}
\end{equation}
Since $0<p_i<1$ and $u \in U_2$, this system has no solutions. Thus, game $(\cG_2,u)$ has no UNE in mixed strategies.
\end{proof}

Similar arguments provide an alternative proof for Proposition \ref{p3_}. 

Let $p =(p_1, p_2)$ be a uniform NE in the game $(\cG_2',u)$  under the Markov realization. Denote the expected payoff   function of player $i$ by $F_{i}$, $i=1, 2$. If $p_1=p_2=1$, the probability of cycling is 1. Otherwise, from \eqref{lim_d2} we obtain
\begin{equation}
\begin{aligned}
    F_{1} = &\frac{q_1[(1-p_1)u_{1}(a_{1})+p_1(1-p_2)u_{1}(a_{2})]+q_2[(1-p_1)p_2u_{1}(a_{1})+(1-p_2)u_{1}(a_{2})]}{1-p_1p_2},\\
    F_{2} = &\frac{q_1[(1-p_1)u_{2}(a_{1})+p_1(1-p_2)u_{2}(a_{2})]+q_2[(1-p_1)p_2u_{2}(a_{1})+(1-p_2)u_{2}(a_{2})]}{1-p_1p_2}.
     \nonumber
\end{aligned}
\end{equation}
Relations \eqref{part_deriv} have the following form in this case:
\begin{equation}
\begin{aligned}
\frac{(q_1 + p_2q_2)(u_{1}(a_{1}) - u_{1}(a_{2}))(1-p_2)}{(1-p_1p_2)^2}=0,\\
\frac{(q_2 + p_1q_1)(u_{2}(a_{1}) - u_{2}(a_{2}))(1-p_1)}{(1-p_1p_2)^2}=0.
\nonumber
\end{aligned}
\end{equation}
Since $u \in U_2$, and for $i = 1, 2$, both $q_i$ cannot be $0$ and by Lemma \ref{l1}, $0 < p_i < 1$, this system has no solutions. 
Thus, $(\cG_2',u)$ has no NE in mixed strategies.
\qed

\begin{proof}[Proof of Proposition \ref{p4}]
Let $p =(p_1, p_2, p_3)$ be a uniform NE in the game $(\cG_3,u)$ under the Markov realization. 
If $p_1=p_2=p_3=1$, the probability of cycling is 1. 
Otherwise, assuming that the initial positions are $v_1$, $v_2$ or $v_3$, 
the limiting distributions are given by \eqref{lim_d3} and 
\begin{equation}
\begin{aligned}
    F_{11} = &\frac{(1-p_1)u_{1}(a_{1})+p_1(1-p_2)u_{1}(a_{2})+p_1p_2(1-p_3)u_{1}(a_{3})}{1-p_1p_2p_3},\\
    F_{21} = &\frac{(1-p_2)u_{1}(a_{2})+p_2(1-p_3)u_{1}(a_{3})+p_2p_3(1-p_1)u_{1}(a_{1})}{1-p_1p_2p_3},\\
    F_{31} = &\frac{(1-p_3)u_{1}(a_{3})+p_3(1-p_1)u_{1}(a_{1})+p_1p_3(1-p_2)u_{1}(a_{2})}{1-p_1p_2p_3},\\
    F_{12} = &\frac{(1-p_1)u_{2}(a_{1})+p_1(1-p_2)u_{2}(a_{2})+p_1p_2(1-p_3)u_{2}(a_{3})}{1-p_1p_2p_3},\\
    F_{22} = &\frac{(1-p_2)u_{2}(a_{2})+p_2(1-p_3)u_{2}(a_{3})+p_2p_3(1-p_1)u_{2}(a_{1})}{1-p_1p_2p_3},\\
    F_{32} = &\frac{(1-p_3)u_{2}(a_{3})+p_3(1-p_1)u_{2}(a_{1})+p_1p_3(1-p_2)u_{2}(a_{2})}{1-p_1p_2p_3},\\
    F_{13} = &\frac{(1-p_1)u_{3}(a_{1})+p_1(1-p_2)u_{3}(a_{2})+p_1p_2(1-p_3)u_{3}(a_{3})}{1-p_1p_2p_3},\\
    F_{23} = &\frac{(1-p_2)u_{3}(a_{2})+p_2(1-p_3)u_{3}(a_{3})+p_2p_3(1-p_1)u_{3}(a_{1})}{1-p_1p_2p_3},\\
    F_{33} = &\frac{(1-p_3)u_{3}(a_{3})+p_3(1-p_1)u_{3}(a_{1})+p_1p_3(1-p_2)u_{3}(a_{2})}{1-p_1p_2p_3}.
   \nonumber
\end{aligned}
\end{equation}
Relations \eqref{part_deriv} have the following form in this case:
\begin{equation}
\begin{aligned}
\frac{-(u_{1}(a_{1}) - u_{1}(a_{2}) + p_2u_{1}(a_{2}) - p_2u_{1}(a_{3}) - p_2p_3u_{1}(a_{1}) + p_2p_3u_{1}(a_{3}))}{(1-p_1p_2p_3)^2}=0,\\
\frac{-(p_2p_3(u_{1}(a_{1}) - u_{1}(a_{2}) + p_2u_{1}(a_{2}) - p_2u_{1}(a_{3}) - p_2p_3u_{1}(a_{1}) + p_2p_3u_{1}(a_{3})))}{(1-p_1p_2p_3)^2}=0,\\
\frac{-(p_3(u_{1}(a_{1}) - u_{1}(a_{2}) + p_2u_{1}(a_{2}) - p_2u_{1}(a_{3}) - p_2p_3u_{1}(a_{1}) + p_2p_3u_{1}(a_{3})))}{(1-p_1p_2p_3)^2}=0,\\
\frac{-(u_{2}(a_{2}) - u_{2}(a_{3}) - p_3u_{2}(a_{1}) + p_3u_{2}(a_{3}) + p_1p_3u_{2}(a_{1}) - p_1p_3u_{2}(a_{2}))}{(1-p_1p_2p_3)^2}=0,\\
\frac{-(p_1p_3(u_{2}(a_{2}) - u_{2}(a_{3}) - p_3u_{2}(a_{1}) + p_3u_{2}(a_{3}) + p_1p_3u_{2}(a_{1}) - p_1p_3u_{2}(a_{2})))}{(1-p_1p_2p_3)^2}=0,\\
\frac{-(p_1(u_{2}(a_{2}) - u_{2}(a_{3}) - p_3u_{2}(a_{1}) + p_3u_{2}(a_{3}) + p_1p_3u_{2}(a_{1}) - p_1p_3u_{2}(a_{2})))}{(1-p_1p_2p_3)^2}=0,\\
\frac{-(u_{3}(a_{3}) - u_{3}(a_{1}) + p_1u_{3}(a_{1}) - p_1u_{3}(a_{2}) - p_1p_2u_{3}(a_{3}) + p_1p_2u_{3}(a_{2}))}{(1-p_1p_2p_3)^2}=0,\\
\frac{-(p_1p_2(u_{3}(a_{3}) - u_{3}(a_{1}) + p_1u_{3}(a_{1}) - p_1u_{3}(a_{2}) + p_1p_2u_{3}(a_{2}) - p_1p_2u_{3}(a_{3})))}{(1-p_1p_2p_3)^2}=0,\\
\frac{-(p_2(u_{3}(a_{3}) - u_{3}(a_{1}) + p_1u_{3}(a_{1}) - p_1u_{3}(a_{2}) - p_1p_2u_{3}(a_{3}) + p_1p_2u_{3}(a_{2})))}{(1-p_1p_2p_3)^2}=0.\nonumber
\end{aligned}
\end{equation}
Since $p_i>0$ and $q_i$ are not all equal to 0 for $i=1, 2, 3$, we have \begin{equation}
\begin{aligned}
u_{1}(a_{1}) - u_{1}(a_{2}) + p_2u_{1}(a_{2}) - p_2u_{1}(a_{3}) - p_2p_3u_{1}(a_{1}) + p_2p_3u_{1}(a_{3})=0,\\
u_{2}(a_{2}) - u_{2}(a_{3}) - p_3u_{2}(a_{1}) + p_3u_{2}(a_{3}) + p_1p_3u_{2}(a_{1}) - p_1p_3u_{2}(a_{2})=0,\\
u_{3}(a_{3}) - u_{3}(a_{1}) + p_1u_{3}(a_{1}) - p_1u_{3}(a_{2}) - p_1p_2u_{3}(a_{3}) + p_1p_2u_{3}(a_{2})=0,
\label{deriv_simpl9}
\end{aligned}
\end{equation}
and using \eqref{mu} we transform equations \eqref{deriv_simpl9} to
\begin{equation}
\begin{aligned}
\mu_1(1-p_2) = p_2(1-p_3), \\
\mu_2(1-p_3) = p_3(1-p_1),\\
\mu_3(1-p_1) = p_1(1-p_2).
\label{deriv_mu}
\end{aligned}
\end{equation}
Recall that $0<p_i<1$ for $i=1, 2, 3$, by Lemma \ref{l1}. Solving \eqref{deriv_mu} WRT $p_i$ yields \eqref{p_i}, provided $p_i>0$, and each equality of \eqref{p_i} implies that $\mu_1 \mu_2 \mu_3 <1$, provided $p_i<1$.
\end{proof}

\begin{proof}[Proof of Proposition \ref{p5_}]
Let $p =(p_1, p_2, p_3)$ be a uniform NE 
in the game $(\cG'_3,u)$ under Markov realization. 
Let the expected payoff of player $i$ be denoted by $F_i$, $i=1, ..., 3$. If $p_1=p_2=p_3=1$, the probability of a cycle is 1. Otherwise, 
\begin{equation}
\begin{aligned}
    F_1 = &\frac{q_1[(1-p_1)u_{1}(a_{1})+p_1(1-p_2)u_{1}(a_{2})+p_1p_2(1-p_3)u_{1}(a_{3})]}{1-p_1p_2p_3}\\
    &+\frac{q_2[(1-p_1)p_2p_3u_{1}(a_{1})+(1-p_2)u_{1}(a_{2})+p_2(1-p_3)u_{1}(a_{3})]}{1-p_1p_2p_3}\\
    &+\frac{q_3[p_3(1-p_1)u_{1}(a_{1})+p_1p_3(1-p_2)u_{1}(a_{2})+(1-p_3)u_{1}(a_{3})]}{1-p_1p_2p_3},\\
    F_2 = &\frac{q_1[(1-p_1)u_{1}(a_{1})+p_1(1-p_2)u_{1}(a_{2})+p_1p_2(1-p_3)u_{1}(a_{3})]}{1-p_1p_2p_3}\\
    &+\frac{q_2[(1-p_1)u_{1}(a_{1})p_2p_3+(1-p_2)u_{1}(a_{2})+p_2(1-p_3)u_{1}(a_{3})]}{1-p_1p_2p_3}\\
    &+\frac{q_3[p_3(1-p_1)u_{1}(a_{1})+p_1p_3(1-p_2)u_{1}(a_{2})+(1-p_3)u_{1}(a_{3})]}{1-p_1p_2p_3},\\
    F_3 = &\frac{q_1[(1-p_1)u_{1}(a_{1})+p_1(1-p_2)u_{1}(a_{2})+p_1p_2(1-p_3)u_{1}(a_{3})]}{1-p_1p_2p_3}\\
    &+\frac{q_2[(1-p_1)u_{1}(a_{1})p_2p_3+(1-p_2)u_{1}(a_{2})+p_2(1-p_3)u_{1}(a_{3})]}{1-p_1p_2p_3}\\
    &+\frac{q_3[p_3(1-p_1)u_{1}(a_{1})+p_1p_3(1-p_2)u_{1}(a_{2})+(1-p_3)u_{1}(a_{3})]}{1-p_1p_2p_3}.\nonumber
\end{aligned}
\end{equation}
Relations \eqref{part_deriv} have the following form in this case:
\begin{equation}
\begin{aligned}
\frac{-(q_1  + q_2p_2p_3+ q_3p_3)(u_{1}(a_{1}) - u_{1}(a_{2}) + p_2u_{1}(a_{2}) - p_2u_{1}(a_{3}) - p_2p_3u_{1}(a_{1}) + p_2p_3u_{1}(a_{3}))}{(1-p_1p_2p_3)^2} = 0,\\
\frac{-(q_2 + q_3p_1p_3 + q_1p_1)(u_{2}(a_{2}) - u_{2}(a_{3}) - p_3u_{2}(a_{1}) + p_3u_{2}(a_{3}) + p_1p_3u_{2}(a_{1}) - p_1p_3u_{2}(a_{2}))}{(1-p_1p_2p_3)^2}=0,\\
\frac{-(q_3 + q_1p_1p_2+ q_2p_2)(u_{3}(a_{1}) - u_{3}(a_{3}) - p_1u_{3}(a_{1}) + p_1u_{3}(a_{2}) - p_1p_2u_{3}(a_{2}) + p_1p_2u_{3}(a_{3}))}{(1-p_1p_2p_3)^2}=0. \label{deriv}
\end{aligned}
\end{equation}
Obviously, for $i=1, 2, 3$, not all $q_i$ are $0$, since $q_1+q_2+q_3=1$, and furthermore, by Lemma \ref{l1}, $0< p_i<1$. Hence, in the LHS of each equation in \eqref{deriv} the denominator and all three first factors are not $0$. Therefore, all three second factors are $0$, which exactly means \eqref{deriv_simpl9}. As before, using \eqref{mu} we transform \eqref{deriv_simpl9} to \eqref{deriv_mu} and obtain \eqref{p_i}, and conditions $p_i<1$ for $i=1, 2, 3$ imply that $\mu_1 \mu_2 \mu_3 <1$. Thus, if $\mu_1\mu_2\mu_3<1$, \eqref{p_i} defines a unique NE, otherwise, if $\mu_1\mu_2\mu_3 \ge 1$, there is no NE.
\end{proof}

Let us note that the above proof works for any distribution $q(v_0)$, not necessarily strictly positive.

\subsection{The a-priori realization}
\label{ss42}

\begin{proof}[Proof of Proposition \ref{p6}]
For game $(\cG_2, u)$ the limiting a-priori distributions for the outcomes $(a_1, a_2, c)$  WRT initial positions $v_1$ and $v_2$, are given by \eqref{lim_d2_ap}. In particular, \eqref{lim_d2_ap} implies that the expected payoffs 
$F_{12}$ and $F_{21}$ are
\begin{equation}
\begin{aligned}
    F_{12} = &(1-p_1)u_{2}(a_{1})+p_1(1-p_2)u_{2}(a_{2})+p_1p_2u_{2}(c),\\
    F_{21} = &(1-p_2)u_{1}(a_{2})+p_2(1-p_1)u_{1}(a_{1})+p_1p_2u_{1}(c).\nonumber
    \label{2pl_payoffs_ap}
\end{aligned}
\end{equation}
Then, the equations \eqref{part_deriv} have the following form:
\begin{equation}
\begin{aligned}
p_1(u_{2}(c)-u_{2}(a_{2}))=0,\\
p_2(u_{1}(c)-u_{1}(a_{1}))=0.\nonumber
\end{aligned}
\end{equation}
If $u \in U_2$, the above system of equations has no solutions, and hence, game  $(\cG_2, u)$ has no UNE in mixed strategies under the a-priori realization.

\medskip

For game $(\cG_3, u)$ the limiting a-priori distributions on the outcomes $(a_1, a_2, a_3, c)$, WRT initial positions $v_1$, $v_2$, and $v_3$, are given by \eqref{lim_d3_ap}.

In particular, the expected mean payoffs 
$F_{21}$, $F_{32}$ and $F_{13}$ are
\begin{equation}
\begin{aligned}
    F_{21} = &(1-p_2)u_{1}(a_{2})+p_2(1-p_3)u_{1}(a_{3})+p_2p_3(1-p_1)u_{1}(a_{1})+p_1p_2p_3u_{1}(c),\\
    F_{32} = &(1-p_3)u_{2}(a_{3})+p_3(1-p_1)u_{2}(a_{1})+p_1p_3(1-p_2)u_{2}(a_{2})+p_1p_2p_3u_{2}(c),\\
    F_{13} = &(1-p_1)u_{3}(a_{1})+p_1(1-p_2)u_{3}(a_{2})+p_1p_2(1-p_3)u_{3}(a_{3})+p_1p_2p_3u_{3}(c).
   \nonumber
\end{aligned}
\end{equation}
The equations \eqref{part_deriv} for them turn into
\begin{equation}
\begin{aligned}
p_2p_3(u_{1}(c)-u_{1}(a_{1}))=0,\\
p_1p_3(u_{2}(c)-u_{2}(a_{2}))=0,\\
p_1p_2(u_{3}(c)-u_{3}(a_{3}))=0.\nonumber
\end{aligned}
\end{equation}
These three equation together contradict Lemma \ref{l1} ($0<p_i<1$), 
when $u \in U_3$. Thus, in this case, 
game $(\cG_3, u)$ has no UNE in mixed strategies under the a-priori realization.
\end{proof}

\section{Interpretation of two main examples} 
\label{s7} 
\subsection{Game $(\cG_2,u)$ with $u \in U_2$} 
Two mechanics $M_1$ and $M_2$ may replace a device in their garage.
There are two options of such replacement: $a_1$  and  $a_2$.
Both prefer  $a_1$  to  $a_2$, so the solution seems obvious.
Yet, there is a third option, $c$:
they do not replace device at all unless they come to consensus.
For  $M_1$  outcome $c$  is the best option: 
better than  $a_1$  (he prefers to save), while
for  $M_2$  $c$  is the worst option: worse than  $a_2$.
They negotiate in pure strategies 
in accordance with the game structure $\cG_2$ on Figure 1.

Suppose $M_1$ makes a move  $(v_1, a_1)$, thus,
agreeing to buy device $a_1$.
Then, naturally,  $M_2$  supports  $M_1$
by making move  $(v_2, v_1)$.
Yet, $M_1$ can improve the obtained situation
$((v_1, a_1),(v_2, v_1))$ for himself rejecting  $a_1$;
that is, he switches  from   $(v_1, a_1)$  to  $(v_1, v_2)$
thus getting  $c$, which is best for him.
This happens for any initial position:
$v_0 = v_1$ or $v_0 =   v_2$.

Recall that  $c$  is the worst outcome for  $M_2$,
so he is unhappy and will improve for himself the current situation
$((v_1, v_2),(v_2, v_1))$ by switching from $(v_2, v_1)$  to $(v_2, a_2)$
and getting  $a_2$  instead of  $c$.
Again, this happens for any initial position:
$v_0 = v_1$ or $v_0 =   v_2$.

Recall that  $a_2$  is the worst outcome  for  $M_1$,
so he is unhappy and will improve for himself the situation
$((v_1, v_2),(v_2, a_2))$ switching
from   $(v_1, v_2)$  to $(v_1, a_1)$  and getting  $a_1$,
at  least when the play begins in  $v_1$.
If it begins in  $v_2$, outcome  $a_2$  remains.
Nevertheless,  $M_1$ makes a strict improvement when $v_0 = v_1$
and he gets the same result when  $v_0 = v_2$.

Finally, $M_2$ can improve the obtained situation
$((v_1, a_1),(v_2, a_2))$  for himself,
switching from  $(v_2, a_2)$  to  $(v_2, v_1)$. 
At  least,  $a_2$  is replaced by  $a_1$
when the play begins in  $v_2$, and
if it begins in  $v_1$  then outcome $a_1$  remains.
Nevertheless,  $M_2$ makes a strict improvement when $v_0 = v_2$
and he gets the same result when  $v_0 = v_2$.


\subsection{Game $(\cG_3,u)$ with $u \in U_3$}  

\subsubsection*{Behavioral interpretation}
Once upon a time there was a family:
grandmother (GM), mother (M), and little girl
(LG, not too little, yet) corresponding to players $1,2,$ and $3$.
The family has a work to do, say, cleaning, washing, or shopping.
Each player can terminate, which means to do the work herself.
This is the second best outcome for each.
Alternatively, each can follow the 3-cycle,
thus, asking the next player to do the work,
in the cyclic order: GM, M, LG.
The  best outcome  for each is
when the next player does the work.
The third best is  when the previous will.
Finally,  $c$  means that nobody did the work, which
is the worst outcome for all.
The following psychological motivation can be suggested.
\begin{itemize}
\item[]
GM  prefers  M  to work, but she pampers LG and
would prefer to work herself instead of her.
\item[]
M prefers LG to work, but has a mercy for GM and
would prefer to replace her.
\item[]
LG, who is already spoiled by GM, prefers her to work,
but not  M, because
in this case  M  may get angry and punish LG somehow in the future.
\end{itemize}

\subsubsection*{Financial interpretation}
Two projects are considered:
\begin{enumerate}[(i)]
\item constructing a bridge 
across the Raritan river in Middlesex County, NJ;
\item including this bridge into a highway (route 18) in future.
\end{enumerate}

Project (ii) is essentially more expensive than  (i).
Only  (i) is under consideration at the present.
Three players are Local (L), State of New Jersey (S), and Federal (F) governments.
All are interested  in projects (i) and (ii),
but also in saving money from their budgets.
Part (ii) is too expensive for L, so either S or F  pays  for it;
S  could  pay for  (i)  or  (ii)  but not for both;
F  could pay for both, but in this case
(ii) will be started only in  8-9 years after  (i).
(You are not alone!)
Otherwise, if  L  or S  pays for (i),
then (ii)  can be started much sooner, say, in 1-2 years.
The big delay is OK with  S, but not with  L.
Both are happy to provide convenient transit,
but the bridge, not included in a highway
will be served only by local roads and 
result in frequent local traffic jams.

\smallskip

Each player can terminate, which means paying for (i),
or refuse to pay, asking the next player to do it; 
in the cyclic order L, S, F.

The best for L  if S pays for (i);
then  F  will pay for  (ii) in 1-2  years.
Yet, if  F  pays for  (i)  then (ii) will be delayed;
so  L  would prefer to pay for (i).

The best for    S  if  F  pays for (i), and then,
in 8-9 years, for (ii).
If  L  pays for (i) then  S  will have
to pay for  (ii) in  1-2 years;
so  S  would rather pay for (i) now, which is much cheaper.

The best for  F  if  L  pays for (i) and then,
in 1-2  years, S   will pay for  (ii).
If  S  pays for (i)  then
F  will have to  pay for  (ii) in  1-2 years;
so  F would rather pay for (i) now and for  (ii) in 8-9  years.

For all three parties $c$  is the worst outcome:
if all refuse to pay then projects (i) and (ii)
will not be realized.

\begin{thebibliography}{99}

\bibitem{AGH10}
D. Andersson, V. Gurvich, and T. D. Hansen,
On acyclicity of games with cycles,
Discrete Applied Mathematics 158:10 (2010) 1049--1063.

\bibitem{Aum64}
R. Aumann,
Mixed and behavior strategies in infinite extensive games,
in M. Dresher, L.S. Shapley, and A.W. Tucker  (eds.),
Advances in Game Theory, Annals of Mathematics Studies, 52,
Princeton, NJ, Princeton University Press (1964) 627--650;
ISBN 9780691079028.

\bibitem{Bla62}
D. Blackwell, Discrete dynamic programming. Ann. Math. Statist. 33 (1962) 719--726.

\bibitem{BEGM12}
E. Boros, K. Elbassioni, V. Gurvich, and K. Makino,
On Nash Equilibria and Improvement Cycles 
in Pure Positional Strategies
for Chess-like and Backgammon-like n-person Games,
Discrete Math. 312:4 (2012) 772--788.

\bibitem{BEGM17} 
E. Boros, K. Elbassioni, V. Gurvich, and K. Makino, 
A nested family of $k$-total effective rewards 
for positional games,
Int. J. Game Theory 46:1 (2017) 263--293.

\bibitem{BFGV22} 
E. Boros, P. G. Franciosa, V. Gurvich, and M. N. Vyalyi, 
Deterministic $n$-person shortest path and 
terminal games on symmetric digraphs have 
Nash equilibria in pure stationary strategies, 
Preprint at arxiv.org/abs/2202.11554 (2022),  
submitted to Int. J. Game Theory.  

\bibitem{BG03}
E. Boros and V. Gurvich,
On Nash-solvability in pure stationary strategies
of finite games with perfect information 
which may have cycles,
Mathematical Social Sciences  46:2 (2003) 207--241.

\bibitem{BG09}
E. Boros and V. Gurvich,
Why chess and backgammon can be solved
in pure positional uniformly optimal strategies,
RUTCOR Research Report 21-2009, Rutgers University.

\bibitem{BGMS11}
E. Boros, V. Gurvich, K. Makino, and W. Shao,
Nash-solvable two-person symmetric cycle game forms,
Discrete Applied Mathematics 159:15 (2011) 1461–-1487.

\bibitem{BGMOV18}
E. Boros, V. Gurvich, M. Milani\v{c}, 
V. Oudalov, and  J. Vi\v{c}i\v{c},
A three-person deterministic graphical game 
without Nash equilibria,
Discrete Applied Math. 243 (2018) 21--38.

\bibitem{BGY13}
E. Boros, V. Gurvich, and E. Yamangil;
Chess-like games may have no uniform
Nash equilibria even in mixed strategies,
Article ID 534875, Hindawi, Game Theory (2013) 1--10.

\bibitem{CJKLS17}
C. S. Calude, S. Jain, B. Khoussainov, W. Li, and F. Stephan, 
Deciding parity games in quasi-polynomial time, 
in H. Hatami, P. McKenzie, and V. King, editors, 
Proceedings of the 49th Annual ACM SIGACT Symposium 
on Theory of Computing, STOC 2017, 
Montreal, QC, Canada, June 19-23 (2017) 252--263.

\bibitem{DM86}
P. Dasgupta and E. Maskin,
The existence of equilibrium 
in discontinuous economic games, 
Review of Economic Studies 53:1 (1986) 1--26.

\bibitem{EF70}
J. Edmonds and D.R. Fulkerson, Bottleneck extrema,
J. Combinatorial Theory 8 (1970) 299--306.

\bibitem{EM79} 
A. Ehrenfeucht and J. Mycielski, 
Positional strategies for mean payoff games, 
Int. J. Game Theory 8 (1979) 109--113.

\bibitem{Gal53}
D. Gale, A theory of N-person games 
with perfect information,
Proceedings of the National Academy of Sciences 
39:6 (1953) 496-–501.

\bibitem{Gil57}
D. Gillette, Stochastic games 
with zero stop probabilities, 
Contributions to the theory of games, 
Annals of Mathematics Studies 39:3 (1957) 179--187. 

\bibitem{GZ04}
H. Gimbert and W. Zielonka, When can you play positionally? 
Mathematical Foundations of Computer Science,  
Lecture Notes in Computer Science  3153 (2004) 686--697. 

\bibitem{Gur73}
V. A. Gurvich, On theory of multistep games, 
USSR Computational Mathematics and Mathematical Physics 
13:6 (1973) 143--161.

\bibitem{Gur75}
V. A. Gurvich,
The solvability of positional games in pure strategies,
USSR Computational Mathematics and Mathematical Physics 15:2 (1975) 74--87.

\bibitem{Gur88}
V. A. Gurvich, 
A stochastic game with complete information and
without equilibrium situations in pure stationary strategies, 
Russian Math. Surveys 43:2 (1988) 171--172. 

\bibitem{Gur89}
V. A. Gurvich,
Equilibrium in pure strategies,
Soviet Mathematics Doklady 38:3 (1989) 597--602.

\bibitem{Gur90}
V. A. Gurvich,
A theorem on the existence of equilibrium situations in pure
stationary strategies for ergodic extensions
of $(2 \times k)$ bimatrix games, 
Russian Math. Surveys 45:4 (1990) 170--172.

\bibitem{Gur17}
V. Gurvich, 
Backward induction in presence of cycles;
Oxford Journal of Logic and Computation 28:7 (2018) 1635--1646. 

\bibitem{Gur17a}
V. Gurvich, 
Generalizing Gale's theorem on backward induction 
and domination of strategies, 
Preprint at  http://arxiv.org/abs/1711.11353 (2017).

\bibitem{Gur21}
V. Gurvich, 
On Nash-solvability of finite $n$-person deterministic 
graphical games, Catch 22, 
Preprint at https://arxiv.org/abs/2111.06278 (2021).   

\bibitem{Gur21a}
V. Gurvich, 
On Nash-solvability of finite n-person shortest path games, 
bi-shortest path conjecture, 
Preprint at  http://arxiv.org/abs/2111.07177  (2021). 

\bibitem{GK18}
V. Gurvich and G. Koshevoy.
Monotone bargaining is Nash-solvable, 
Discrete Applied Mathematics 250 (2018) 1--15. 

\bibitem{GKK88}
V. Gurvich, A. V. Karzanov, and L. Khachiyan,
Cyclic games and an algorithm to find
minimax cycle means in directed graphs,
USSR Comput. Math. and Math. Phys. 28:5 (1990) 85--91.

\bibitem{GN22}
V. Gurvich and M. Naumova, 
Lexicographically maximal edges of dual hypergraphs and
Nash-solvability of tight game forms, 
Annals of Mathematics and Artificial Intelligence 
2022, available at https://doi.org/10.1007/s10472-022-09820-3, 9 pages. 

\bibitem{GN22a} 
V. Gurvich and M. Naumova, 
Polynomial algorithms computing two lexicographically safe 
Nash equilibria in finite two-person games 
with tight game forms given by oracles; 
Preprint at  https://arxiv.org/abs/2108.05469 , 
published online 25 January 2022, 26 pages. 

\bibitem{GN22b} 
V. Gurvich and M. Naumova, 
On Nash-solvability of finite two-person tight 
vector game forms, 
Preprint at  https://arxiv.org/abs/2204.10241 , 
published online 22 April 2022, 17 pages. 

\bibitem{GO14}
V. Gurvich and V. Oudalov,
On Nash-solvability in pure stationary strategies
of the deterministic n-person games
with perfect information and mean or total effective cost,
Discrete Appl. Math. 167 (2014) 131--143.

\bibitem{GO15}
V. Gurvich and V. Oudalov,
A four-person chess-like game
without Nash equilibria in pure stationary strategies, 
Business Informatics 1:31 (2015) 68--76.

\bibitem{Hov60}
R. A. Howard, Dynamic Programming and Markov Processes, 
The M.I.T. Press, 1960.

\bibitem{HKW83}
A. Hordijk, O. J. Vrieze, and G. L. Wanrooij, 
Semi-markov strategies in stochastic games,
International Journal of Game Theory 12 (1983) 81-–89.

\bibitem{KS60}
J. G. Kemeny and J. L. Snell,
Finite Markov Chains, Springer, 1960.

\bibitem{Kuh50}
H. Kuhn, Extensive games,
Proc. Nat. Acad. Sci. 36 (1950) 286--295.

\bibitem{Kuh53}
H. Kuhn, Extensive games and the problems of information,
Annals Math. Studies 28 (1953) 193--216.

\bibitem{KFSV09}
J. Kuipers, J. Flesch, G. Schoenmakers, and K. Vrieze,
Pure subgame-perfect equilibria in free transition games,
European J. Oper. Res, 199:2 (2009) 442--447.

\bibitem{LL69}
T.M. Liggett and S.A. Lippman, Stochastic games with perfect information 
and time average payoff, SIAM Rev. 11 (1969) 604--607.

\bibitem{MO70}
H. Mine and S. Osaki, Markovian Decision Process,
American Elsevier, New York, NY, USA, 1970.

\bibitem{Mou76}
H. Moulin,  
Prolongement des jeux a 
deux joueurs de somme nulle, 
Une theorie abstraite des duels, 
Memoires de la Societe Mathematique de France,
45 (1976) 5--111;  doi:10.24033/msmf.180

\bibitem{Nas50}
J. Nash,
Equilibrium points in n-person games,
Proc. Nat. Acad. Sci. 36:1 (1950) 48--49.

\bibitem{Nas51}
J. Nash, Non-cooperative games,
Annals of Math. 54:2 (1951) 286--295.

\bibitem{Ros65}
J.B. Rosen,
Existence and uniqueness of equilibrium points 
for concave N-person games,
Econometrica, 33:3 (1965) 520--534.

\bibitem{Sha53}
L. Shapley, Stochastic games, 
Proc. Nat. Acad. Sci. USA 39 (1953) 1095--1100. 

\bibitem{TR97}
F. Thuijsman and  E. S. Raghavan, 
Perfect Information Stochastic Games and Related Classes, 
Int. J.  Game Theory  26:3 (1997) 403--408.

\bibitem{TV87}
F. Thuijsman and O.J.  Vrieze, 
The bad match, a total reward stochastic game. 
Oper. Res. Spektrum 9 (1987) 93--99.

\bibitem{TV98} 
F. Thuijsman and O.J. Vrieze, 
Total reward stochastic games and sensitive average reward strategies. 
J. Optim. Theory Appl. 98 (1998) 175--196.

\bibitem{Was90}
A. Washburn,
Deterministic graphical games, J. Math. Analysis and Appl. 153:1 (1990) 84--96.

\end{thebibliography}
\end{document}